\documentclass[11pt,a4paper]{article}
\usepackage{amsfonts}
\usepackage{amssymb}
\usepackage{graphicx}
\usepackage{amsmath}
\usepackage{verbatim}
\usepackage{color}
\usepackage{pifont}
\usepackage{float}
\usepackage{hyperref}
\usepackage[dvipsnames]{xcolor}
\usepackage[margin=15mm,font=small,labelfont=bf]{caption}
\usepackage[OT2,T2A,T1]{fontenc}
\usepackage[utf8]{inputenc}
\usepackage[russian,french,english]{babel}
\usepackage[toc,title,page]{appendix}

\newtheorem{theorem}{Theorem}

\newtheorem{definition}[theorem]{Definition}

\newtheorem{lemma}[theorem]{Lemma}
\newtheorem{notation}[theorem]{Notation}

\newtheorem{remark}{Remark}[theorem]

\newenvironment{proof}[1][Proof]{\noindent\textbf{#1.} }{\mbox{ } \hfill \rule{0.5em}{0.5em}\medskip}
\numberwithin{equation}{section}

\begin{document}

\title{On a formula for all sets of constant width in 3d}
\author{Bernd Kawohl\footnotemark ,\setcounter{footnote}{0} \ Guido Sweers%
\thanks{%
Dept.~Mathematik \& Informatik, Universit\"{a}t zu K\"{o}ln,
Albertus-Magnus-Platz, 50923 K\"{o}ln, Germany}}
\maketitle

\footnotetext{%
Email: \texttt{kawohl@math.uni-koeln.de gsweers@math.uni-koeln.de}}
\footnotetext{%
Orcid-Id: Kawohl 0000-0003-2918-7318; Sweers 0000-0003-0180-5890\medskip}

\abstract{In the recent paper \lq\lq On a formula for sets of constant width in 2D\rq\rq,
Comm.~Pure Appl.~Anal. 18 (2019), 2117--2131,
we gave a constructive formula for all 2d sets of constant width. Based on this result we derive here
a formula for the parametrization of the boundary of bodies of constant width in
3 dimensions, depending on one function defined on $\mathbb{S}^2$. Each such function gives a minimal
value $r_0$ and for all $r\ge r_0$ one finds a body of constant width $2r$.
Moreover, we show that all bodies of constant width in 3d have such a parametrization.
The last result needs a tool that we describe as \lq shadow domain\rq\ and that is explained in an appendix.
The construction is explicit and and offers a parametrization different from the one given by T. Bayen, T. Lachand-Robert
and \'{E}. Oudet, \lq\lq Analytic parametrization of three-dimensional bodies of constant width\rq\rq\
in Arch.~Ration.~Mech.~Anal., 186 (2007), 225--249.}%
\bigskip

\smallskip


\footnotetext{\textbf{Acknowledgement:} The authors thank Prof.~Hansj\"org
Geiges for pointing out reference \cite{MP} and Ameziane Oumohand M.Sc.~for
\cite{Gr}. Thanks also go to a referee, since the final version benefitted
from the careful and detailed report. \smallskip

It should be mentioned that a first version of the manuscript was essentially completed while
the first author participated in the program \lq\lq Geometric Aspects of
Nonlinear Partial Differential Equations\rq\rq, which was supported by the
Swedish Research Council, at Institut Mittag-Leffler in Djursholm, Sweden
during October 2022.}

\section{Introduction}

For a compact set $G\subset \mathbb{R}^{n}$ one defines its directional
width in direction $\omega \in \mathbb{S}^{n-1}:=\left\{ x\in \mathbb{R}%
^{n};\left\vert x\right\vert =1\right\} $ by%
\begin{equation*}
\boldsymbol{d}_{G}\left( \omega \right) =\max \left\{ \left\langle \omega
,x\right\rangle ;x\in G\right\} -\min \left\{ \left\langle \omega
,x\right\rangle ;x\in G\right\} ,
\end{equation*}%
with $\langle \cdot ,\cdot \rangle $ denoting the standard inner product. If $G$ is
convex and $\boldsymbol{d}_{G}\left( \omega \right) =\boldsymbol{d}_{G}$ is
constant, then $G$ is called a set of constant width. In 3 dimensions a set of constant width
is also called a body of constant width.\medskip

The interest in the subject started with Leonhard Euler, who around 1774
considered 2d curves of constant width, which he called `curva orbiformis'.
He not only studied such sets for 2 dimensions but also gave a formula
describing such curves. See \S 10 of \cite{Eu}. In 3 dimensions a ball is
obviously the classical example of a body of constant width but the famous
Meissner bodies also have this property. See \cite{Me1, Me2} or \cite{KW}.
Quite simple examples can also be constructed by taking a reflection symmetric 2d set
of constant width and rotating it around its line of symmetry. \medskip

Famous mathematicians such as Minkowski \cite{Mi} and Hilbert \cite{HC} were
intrigued by the subject. The first interest of most scholars focused on
deriving properties of such domains. A wonderful survey on sets of constant
width (up to 1983) was provided by Chakerian and Groemer in \cite{CG}, and a
more recent updated and thorough treatment can be found in the book by
Martini, Montejano and Oliveros \cite{MMO}. Let us recall that the 3d
question, motivated by Blaschke's 2d result \cite{Bl}, as to which body of
constant fixed width has the smallest volume or, equivalently, the smallest
surface area, is still open. We will not solve that problem, but will give
an alternative formula for constructing bodies of constant width that might
help. \medskip

Let us recall some known facts about sets of constant width.
Sets of constant width $G$ in $\mathbb{R}^{n}$ are strictly convex and hence
any tangential plane touches $G$ in at most one point. Moreover, for every
boundary point $\boldsymbol{X}$ of $G$ and $r>\boldsymbol{d}_G$ there exists even a ball $B_r$ of
radius $r$ such that $G\subset B_r$ and $\partial G\cap \partial B_r=\boldsymbol{X}$.
So $G$ slides freely in $B_r$, according to the definition of \cite[page 244]{Ho}.

Although the
Gauss-map $\partial \Omega \rightarrow \mathbb{S}^{n-1}$ (the outward normal
on smooth parts of $\partial G$) will not be uniquely defined on edges or corners, the
`inverse' $\gamma _{G}:\mathbb{S}^{n-1}\rightarrow \partial G$ is
well-defined for a strictly convex $G$ and parametrizes $\partial G$. See \cite%
{LO}. In \cite[Theorem 11.1.1]{MMO} one finds, when $G$ is a set of constant
width and $\gamma _{G}\in C^{1}\left( \mathbb{S}^{n-1}\right) $, that
\begin{equation}
\gamma _{G}(\omega )=P_{G}(\omega )\,\omega +\nabla _{P_{G}}(\omega ) \qquad
\text{ for all }\omega \in \mathbb{S}^{n-1},
\label{InvG}
\end{equation}%
where $P_{G}(\omega ):=\max \left\{ \left\langle \omega ,x\right\rangle
;x\in G\right\} $ is the support function and $\nabla _{P_{G}}$ its gradient
along $\mathbb{S}^{n-1}$, i.e. $\left\langle \nabla _{P_{G}}(\omega
),v\right\rangle =\left( dP_{G}(\omega )\right) (v)$ for all $\omega \in
\mathbb{S}^{n-1}$ and $v\in T_{\omega }\mathbb{S}^{n-1}$. Howard in \cite[Corollary 2.6]{Ho}
states that for a body of constant width the support function $P_G$ is of class $C^{1,1}$
and hence $\gamma_G$ is Lipschitz continuous. A direct proof of this Lipschitz-continuity
also follows from Lemma \ref{techlip} below.
Necessary for a set $G$ of constant width $\boldsymbol{d}_{G}$ is that
\begin{equation}
\gamma _{G}(\omega )-\gamma _{G}(-\omega )=\boldsymbol{d}_{G}\ \omega\qquad
\text{ for all }\omega \in \mathbb{S}^{n-1}.
\label{NSCo}
\end{equation}%
Therefore $P_{G}(\omega )+P_{G}(-\omega )=\boldsymbol{d}_{G}$
and $\nabla _{P_{G}}(\omega )=\nabla _{P_{G}}(-\omega )$ for all $\omega \in
\mathbb{S}^{n-1}$.\medskip

In \cite{BLO} a parametrization of sets of constant width is given by using
the so-called median surface, which is parametrized by%
\begin{equation}
M_{G}(\omega ):=\gamma _{G}(\omega )-\tfrac{1}{2}\boldsymbol{d}_{G}\ \omega \qquad
\text{ for all }\omega \in \mathbb{S}^{n-1}.\medskip   \label{Mgam}
\end{equation}%
Writing $x\cdot Y:=\sum_{i}x_{i}Y_{i}$, which may coincide with but will not
be restricted just to the inner product $\left\langle\cdot,\cdot\right\rangle$,
the convexity of $G$ leads to%
\begin{equation}
\left( M_{G}(\hat{\omega})-M_{G}(\omega )\right) \cdot \omega \leq \tfrac{1}{%
4}\boldsymbol{d}_{G}^{2}\,\left\vert \hat{\omega}-\omega \right\vert ^{2}\text{
for all }\omega ,\hat{\omega}\in \mathbb{S}^{n-1},  \label{Mco1}
\end{equation}%
while (\ref{NSCo}) implies%
\begin{equation}
M_{G}(\omega )=M_{G}(-\omega )\text{ for all }\omega \in \mathbb{S}^{n-1}.
\label{Mco2}
\end{equation}

The reverse question would be: can one give criteria on a continuous
function $\gamma :\mathbb{S}^{n-1}\rightarrow \mathbb{R}^{n}$ such that $%
\gamma \left( \mathbb{S}^{n-1}\right) $ parametrizes the boundary of a set
of constant width? An answer is given by Theorem 2 of \cite{BLO}, where it
is stated that for any continuous map $M:\mathbb{S}^{n-1}\rightarrow \mathbb{%
R}^{n}$ and $\alpha >0$, which satisfy%
\begin{equation}
\begin{array}{cl}
M(\omega )=M(-\omega )\text{ } & \text{for all }\omega \in \mathbb{S}^{n-1},
\\
\left( M(\hat{\omega})-M(\omega )\right) \cdot \omega \leq \frac{1}{4}\alpha
^{2}\ \left\vert \hat{\omega}-\omega \right\vert ^{2} & \text{for all }%
\omega ,\hat{\omega}\in \mathbb{S}^{n-1},%
\end{array}
\label{Mco}
\end{equation}%
the set%
\begin{equation}
G:=\left\{ M(\omega )+t\omega \ ;\ \omega \in \mathbb{S}^{n-1},0\leq t\leq
\tfrac{1}{2}\alpha \right\}\label{GfromM}
\end{equation}%
is of constant width $\boldsymbol{d}_{G}:=\alpha $ and $M_{G}(\omega ):=M(\omega
)$. The $\gamma _{G}$ and $M_{G}$ are as in (\ref{Mgam}).

One finds by continuity of $M$ that
\begin{equation}
\partial G\subseteq \left\{ M(\omega )+\tfrac{1}{2}\alpha \ \omega \ ;\ \omega
\in \mathbb{S}^{n-1}\right\}   \label{bopa}
\end{equation}%
and even that the identity holds in (\ref{bopa}).\medskip

Continuity or even differentiability of $M$ by itself is not enough for an $\alpha $
to exist for which (\ref{Mco}) holds. The second condition in (\ref{Mco})
implies the convexity of $G$ from (\ref{GfromM}) and as such it gives a monotonicity for directional
derivatives, hence a necessary one-sided estimate for second derivatives of $%
M$, whenever these exist. In two dimensions, see \cite{KS}, a few simple
conditions on a function in $L^{\infty }( 0,\pi ) $ are necessary
and sufficient in order to have a curve of constant width. The construction
in 2 dimensions is also helpful in 3 dimensions. It will allow us to give a
more explicit formula for all bodies of constant width, which is what we
want to show here.

\section{Two dimensions}

In the last century Hammer and Sobczyk described a construction for 2
dimensions in \cite{HS1, HS2, Ha1}, based on a characterization of what they
called `outwardly simple line families'. More recently a direct concise
formula was given in \cite{KS} to describe all those sets in two dimensions
starting from any $L^{\infty }( 0,\pi ) $-function satisfying 2
equations, namely the ones in (\ref{a2-2d}). Let us recall the 2d formula
from \cite{KS}:

\begin{theorem}[{\protect\cite[Theorem 3.2]{KS}}]
\label{recipe2d} Let $\boldsymbol{x}_{0}\in \mathbb{R}^{2}$, $r\in \mathbb{R}
$ and $a\in L^{\infty }\left( \mathbb{R}\right) $ satisfy%
\begin{align}
& r\geq \left\Vert a\right\Vert _{\infty },  \label{r-2d} \\[0.04in]
& a\left( \varphi +\pi \right) =-a\left( \varphi \right) \text{ for all }%
\varphi ,  \label{a1-2d} \\
\int_{0}^{\pi }& a\left( s\right) \binom{-\sin s}{\cos s}\ ds=\binom{0}{0}.
\label{a2-2d}
\end{align}%
Define the closed curve $\boldsymbol{x}:\left[ 0,2\pi \right] \rightarrow
\mathbb{R}^{2}$ by
\begin{equation}
\boldsymbol{x}\left( \varphi \right) =\boldsymbol{x}_{0}+\int_{0}^{\varphi
}\left( r-a\left( s\right) \right) \binom{-\sin s}{\cos s}ds.  \label{2df}
\end{equation}%
Then $\boldsymbol{x}$ describes the boundary of a set of constant width $2r$.
\end{theorem}

For a simple statement in Theorem \ref{recipe2d} the function $a\in
L^{\infty }\left( 0,\pi \right) $ is extended to $\mathbb{R}$ and such that (%
\ref{r-2d}) and (\ref{a1-2d}) are satisfied. The formula in (\ref{2df})
shows that $\boldsymbol{x}\in C^{0,1}\left( \mathbb{R}\right) $, which is
optimal for $r=\left\Vert a\right\Vert _{\infty }$. For $r>\left\Vert
a\right\Vert _{\infty }$ when considering the set $\boldsymbol{x}(\left[
0,2\pi \right])$ as a curve one finds that $\boldsymbol{x}(\left[ 0,2\pi %
\right]) \in C^{1,1}$. The formula in (\ref{2df}) describes the boundary of
all 2d domains of constant width:

\begin{theorem}[{\protect\cite[Theorem 4.1]{KS}}]\label{Th2}
If $G\subset \mathbb{R}^{2}$ is a closed convex set of constant width $2r$,
then there exists $\boldsymbol{x}_{0}$ and $a$ as in Theorem \ref{recipe2d},
such that $\partial G=\boldsymbol{x}\left( \left[ 0,2\pi \right] \right) $
with $\boldsymbol{x}$ as in (\ref{2df}).
\end{theorem}

The geometric interpretation of the formula in (\ref{2df}) is that $%
\boldsymbol{x}(\varphi)$ and $\boldsymbol{x}(\varphi+\pi)$ describe the ends
of a rotating stick of length $2r$ with the varying point of rotation lying
on the stick by (\ref{r-2d}) and determined by $a(\varphi)$. For these ends
to coincide for $\varphi\in[0,\pi]$ with those for $\varphi\in[\pi,2\pi]$
one needs condition (\ref{a1-2d}). The two equalities in condition $(\ref%
{a2-2d})$ make it a closed curve.

\section{A formula in three dimensions}

There have been previous attempts to provide an explicit construction of all
3d bodies of constant width. In \cite{LO} Lachand-Robert and Oudet present a
geometric construction that generates 3d bodies of constant width from 2d
sets of constant width. This construction, however, does not capture all 3d
bodies of constant width because a counterexample is provided in the paper
\cite{Da} by Danzer, who constructs a body of constant width $\boldsymbol{d}$,
none of whose planar cross-sections have constant width $\boldsymbol{d}$ in two dimensions.
In \cite{MR} Montejano and Roldan-Pensado generalize
the construction of Meissner bodies to generate so-called Meissner
polyhedra. This construction does not generate all 3d bodies either, because
the rotated Reuleaux triangle is a counterexample. As already mentioned Bayen,
Lachand-Robert and Oudet give a description of (all) $n$-dimensional sets of constant
width in \cite[Theorem 2]{BLO}, but the function $M$ has to satisfy a
condition at each point of $\mathbb{S}^{n-1}$. We provide an alternative
construction, based on the method from \cite{KS}, which gives a simpler
condition although more involved than some integral conditions and an $L^\infty$ bound. Indeed, some simple
conditions as in 2d do not seem possible, but our conditions in 3d will come close. \medskip

Our approach uses spherical coordinates in $\mathbb{R}^3$. Indeed, for each fixed angle
$\theta$ we apply the 2d-approach to get a curve parametrized by $\varphi$ of constant width $%
2r$. So as a first step the function $a$ from Theorem \ref{recipe2d} is now depending on $\theta$
\begin{equation}  \label{afunc}
\varphi\mapsto a(\varphi,\theta)\text{ for each }\theta,
\end{equation}
and is used to define a curve $\varphi\mapsto\boldsymbol{x}(\varphi;\theta)$, with $\boldsymbol{x}$ as
in (\ref{2df}) and $\theta$ as a parameter, in the $\theta$-dependent plane
\begin{equation*}
\boldsymbol{X}_0 + \boldsymbol{x}_1(\varphi;\theta)\left(\begin{array}{c}
                               \cos\theta \\
                               \sin\theta \\
                               0
                             \end{array}\right)
                  +\boldsymbol{x}_2(\varphi;\theta) \left(\begin{array}{c}
                               0 \\
                               0\\
                               1
                             \end{array}\right).
\end{equation*}
This first step however does not yet generate a body of constant width.
Whenever $\left\|\partial^2_\theta a\right\|_\infty$ is bounded and when $r$
is large enough, the second step is to apply a unique shift in the perpendicular
$\left(-\sin\theta,\cos\theta,0\right)^\top$-direction for the collection of these rotating 2d-curves.
For the magnitude of the shift we will use $h(\varphi,\theta)$.
The combined result of these two steps will yield a 3d-body of
constant width. Moreover we will show, that not only the result is a body of
constant width but also that each such body can be written this
way.\medskip

Aside from our results from \cite{KS} for two dimensions we will use a
result by Hadwiger in \cite{Ha}, which can be roughly described as follows: convex
bodies in $\mathbb{R}^{n}$ are uniquely determined by the projections in $%
\mathbb{R}^{n-1}$ perpendicular to one fixed direction. The result holds for
$n\geq 4$ and, whenever the one fixed direction is regular, also for $n=3$.
This last addendum is due to \cite{Gr}. Regular means here, that the planes
perpendicular to that fixed direction which touch the convex domain, do that
in precisely one point. Since sets of constant width are necessarily
strictly convex, this is obviously the case for those sets and any choice of
the fixed direction.\medskip

\begin{figure}[h]
\centering
\includegraphics[width=58mm]{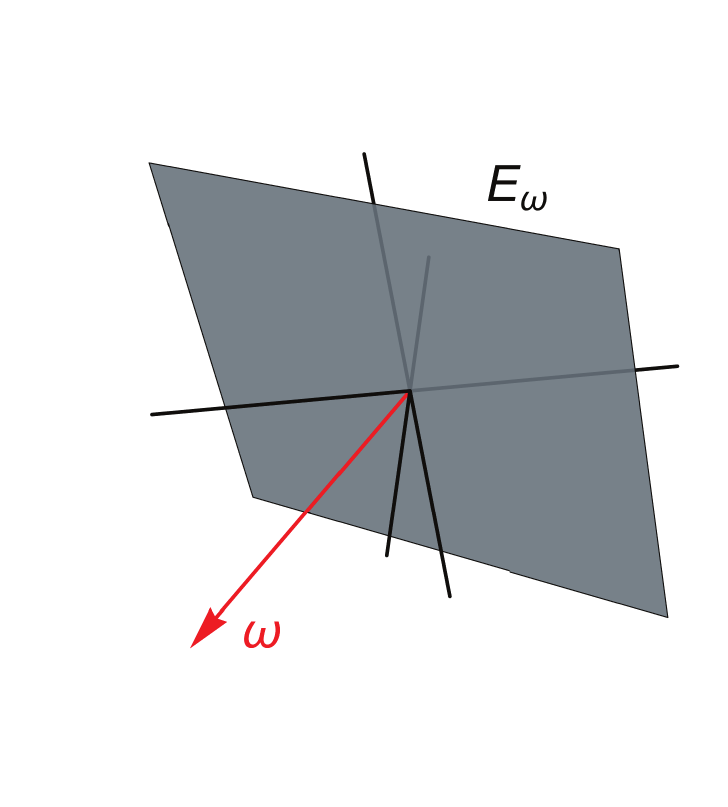}\hspace{9mm}%
\includegraphics[width=6cm]{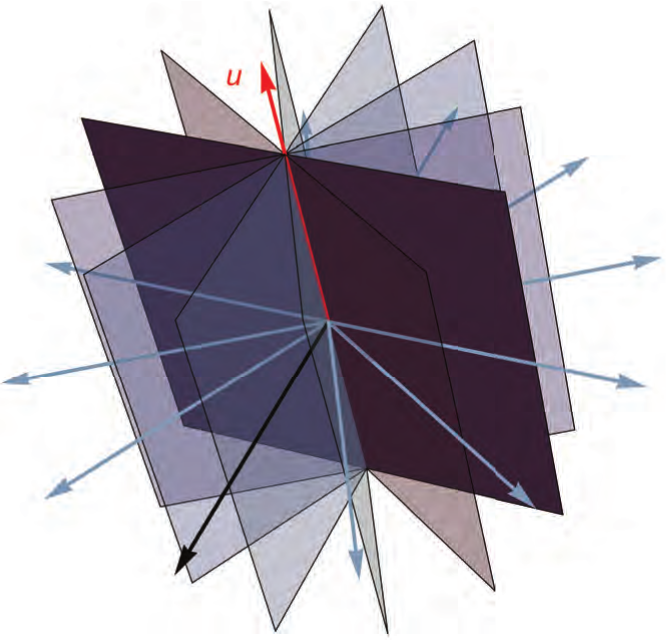}
\caption{The plane $E_{\protect\omega }$ for one $\protect\omega $ and \lq
all\rq\ planes $E_{\protect\omega }$ with $\protect\omega $ such that $%
\left\langle \protect\omega ,u\right\rangle =0$. Those $E_{\protect\omega }$
contain $u$ as a common direction.}
\label{waaier}
\end{figure}

Let us define for $\omega \in \mathbb{S}^{2}$ the orthogonal projection $%
P_{\omega }$ on the plane $E_{\omega }:=\left\{x\in \mathbb{R}^{3};\langle
x,\omega\rangle=0 \right\}$. To exploit the result of Hadwiger we will use
for a fixed $u\in \mathbb{S}^{2}$ all projections in the directions $\omega
\in \mathbb{S}^{2}$ with $\left\langle \omega ,u\right\rangle =0$. See
Fig.~\ref{waaier}. For those $\omega$ we have%
\begin{equation}
P_{\omega }x=\left\langle u,x\right\rangle u+\left\langle u\times \omega
,x\right\rangle \left( u\times \omega \right) .  \label{Pom}
\end{equation}%
For later use we need to identify the projections on $E_{\omega }$ with
coordinates in $\mathbb{R}^{2}$ through%
\begin{equation}
\widehat{P}_{\omega }x=\left(
\begin{array}{c}
\left\langle u,x\right\rangle \\
\left\langle u\times \omega ,x\right\rangle%
\end{array}%
\right) .  \label{Poms}
\end{equation}

We may now explain the result by Hadwiger in \cite{Ha} in more detail. He
proved that for two convex bodies $G_{1}$ and $G_{2}$ in $\mathbb{R}^{3}$
the following holds.

\begin{itemize}
\item \textit{If }$P_{u}G_{1}\simeq P_{u}G_{2}$ \textit{and} $P_{\omega
}G_{1}\simeq P_{\omega }G_{2}$ \textit{for all} $\omega \in \mathbb{S}^{2}$
\textit{with} $\left\langle \omega ,u\right\rangle =0$\textit{, then }$%
G_{1}\simeq G_{2}$.
\end{itemize}

Here $A\simeq B$ means that $A$ equals $B$ after a translation. In other
words, there is a fixed $v\in \mathbb{R}^{3}$ such that $A=v+B$. Groemer
showed in \cite{Gr} that one could drop the condition $P_{u}G_{1}\simeq
P_{u}G_{2}$, whenever $u$ is a regular direction for $G_{1}$. Here regular
means that $\max \left\{ \left\langle u,x\right\rangle ;x\in G_{1}\right\} $
is attained for a unique $x\in G_{1}$. Since domains $G$ of constant width
are precisely those domains for which
\begin{equation*}
G^{\ast }:=\tfrac{1}{2}G+\tfrac{1}{2}\left( -G\right) :=\left\{ \tfrac{1}{2}%
x-\tfrac{1}{2}y;x,y\in G\right\}
\end{equation*}%
is a ball, which has only regular directions, one finds that $\big(\widehat{P%
}_{\omega }G\big)^{\ast }$ is a disc for all $\omega \in \mathbb{S}^{2}$
with $\left\langle \omega ,u\right\rangle =0$, if and only if $G^{\ast }$ is
a ball. Necessarily those discs and the ball have the same radius. This
implies that a convex closed set\ $G\subset \mathbb{R}^{3}$ is a body of
constant width if and only if there is a direction $u\in \mathbb{S}^{2}$,
such that for some fixed $\rho >0$ one finds
\begin{equation*}
\big(\widehat{P}_{\omega }G\big)^{\ast }\simeq D_{\rho }:=\left\{ y\in
\mathbb{R}^{2};\left\vert y\right\vert \leq \rho \right\} \text{ for all }%
\omega \in \mathbb{S}^{2}\text{ with }\left\langle \omega ,u\right\rangle =0.
\end{equation*}%
This means that all those $P_{\omega }G$ should be two-dimensional convex
sets of constant width $\rho $. So by taking $u=\left( 1,0,0\right) $ we
find that the boundary of $P_{\omega }G$ is described by (\ref{2df}) with
some $a$ depending on $\omega $. This leads us to the result in Theorem \ref%
{recipe3d} 

\begin{notation}
\label{param}We parametrize $\mathbb{S}^{2}=\boldsymbol{U}(\mathbb{R}^{2})$
by
\begin{equation}
\omega = \boldsymbol{U}\left( \varphi ,\theta \right) :=\left(
\begin{array}{c}
\sin \varphi \cos \theta \\
\sin \varphi \sin \theta \\
\cos \varphi%
\end{array}%
\right) .  \label{defpar}
\end{equation}
\end{notation}

This is the standard parametrization with $\varphi $ the angle between $\omega$ and the
positive $z$-axis and $\theta $ the counterclockwise angle of the projection
on the $xy$-plane with the $x$-axis, viewed from the positive $z$-axis.
Obviously this parametrization is not unique as we may restrict $%
(\varphi,\theta)$ to some subset of $\mathbb{R}^2$.\medskip

We may define a convenient $\varphi ,\theta $-dependent orthonormal basis,
first for $\sin \varphi\ne 0$,
\begin{align}
\left\{ \boldsymbol{U}\left( \varphi ,\theta \right) ,\boldsymbol{U}%
_{\varphi }\left( \varphi ,\theta \right) ,\frac{\boldsymbol{U}_{\theta
}\left( \varphi ,\theta \right) }{\sin \varphi }\right\}&= \left\{ \left(
\begin{array}{c}
\sin \varphi \cos \theta \\
\sin \varphi \sin \theta \\
\cos \varphi%
\end{array}%
\right) ,\left(
\begin{array}{c}
\cos \varphi \cos \theta \\
\cos \varphi \sin \theta \\
-\sin \varphi%
\end{array}%
\right) ,\left(
\begin{array}{c}
-\sin \theta \\
\cos \theta \\
0%
\end{array}%
\right) \right\}  \notag \\
&=: \left\{ \boldsymbol{U}\left( \varphi ,\theta \right) ,\boldsymbol{V}%
\left( \varphi ,\theta \right) ,\boldsymbol{W}\left( \theta \right) \right\},
\label{movbas}
\end{align}%
with the expression in the middle showing the obvious extension when $\sin
\varphi =0$.\medskip

Any function $(\varphi,\theta)\mapsto v(\varphi,\theta):\mathbb{R}^2\to%
\mathbb{R}$ that is used to define a quantity on $\mathbb{S}^2$ necessarily
has to possess the obvious periodicity properties as well as some
compatibility conditions. The relations, which the function $a$ from (\ref%
{afunc}) has to satisfy, are more subtle. For $\varphi\not\in\left\{0,\pi%
\right\}$ the value $r-a(\varphi,\theta)$ coincides with the inverse
curvature in the $\varphi$-direction. There is however a peculiarity at the
north- and southpole, where the curvature in any(!)~direction is given by $%
\left(r-a(0,\theta)\right)^{-1}$, respectively $\left(r-a(\pi,\theta)%
\right)^{-1}$, through varying $\theta$. This leads to the following
definition with a distinction between pure periodicity and what we call
compatibility, both derived from $\boldsymbol{U}(\varphi,\theta)=\boldsymbol{U}(\hat{\varphi},\hat{\theta})$:

\begin{definition}
\label{funnydifdef}For a function $f:\mathbb{R}^2\rightarrow \mathbb{R}$ we
say that:

\begin{itemize}
\item $f$ satisfies the periodicity conditions for $\mathbb{S}^2$, if
\begin{align}
f(\hat{\varphi},\hat{\theta})&=f(\varphi ,\theta )\ \text{ for all }\ \hat{%
\varphi}-\varphi ,\hat{\theta}-\theta \in 2\pi \mathbb{Z}\ \text{ and}  \label{perio1} \\
f(\hat{\varphi},\hat{\theta})&=f(\varphi ,\theta )\ \text{ for all }\ \hat{%
\varphi}+\varphi ,\hat{\theta}-\theta +\pi \in 2\pi \mathbb{Z};
\label{perio2}
\end{align}

\item $f$ satisfies the compatibility conditions for the poles of $\mathbb{S}%
^2$, if
\begin{equation}
f(0,\theta )=f(0,0)\ \text{ and }\ f(\pi ,\theta )=f(\pi ,0)\ \text{ for all
}\ \theta \in \mathbb{R}.  \label{compa}
\end{equation}
\end{itemize}

\noindent Suppose that $B(\mathbb{R}^2)$ is some function space. We write:

\begin{itemize}
\item $f\in B_{\mathrm{p}}(\mathbb{R}^2)$, whenever $f\in B(\mathbb{R}^2)$
and satisfies (\ref{perio1}) and (\ref{perio2});

\item $f\in B_{\mathrm{p,c}}(\mathbb{R}^2)$, whenever $f\in B(\mathbb{R}^2)$
and satisfies (\ref{perio1}), (\ref{perio2}) and (\ref{compa}).
\end{itemize}
With obvious changes we use the similar notations for a vector-valued $\boldsymbol{F}:\mathbb{R}^2\rightarrow \mathbb{R}^3$ whenever $f=\boldsymbol{F}_i$ satisfies the required properties for all $i\in\{ 1,2,3\}$.
\end{definition}

One usually restricts $\mathbb{R}^{2}$ to $\left[ 0,\pi \right] \times \left[
0,2\pi \right] $ to have a unique parametrization at least for the interior
points and with some compatibility assumptions at its boundary, but here it
will be more convenient to take
\begin{equation}\label{es}
\mathsf{S}=\left[ 0,2\pi \right] \times \left[ 0,\pi \right].
\end{equation}

As in the 2d-case the function $a$ from (\ref{afunc}) that we use is such
that at opposite points of $\mathbb{S}^{2}$ the value is opposite: $a\left(
\varphi ,\theta \right) =-a\left( \varphi +\pi ,\theta \right) =-a\left( \pi
-\varphi ,\theta +\pi \right) $. Hence $a$ is completely defined by its
values on $\left[ 0,\pi \right) \times \left[ 0,\pi \right)
$.

\begin{theorem}[Constructing bodies of constant width]
\label{recipe3d} Suppose that $a\in C_\mathrm{p}^2(\mathbb{R}^2)$ satisfies
\begin{align}
a\left( \varphi ,\theta \right) =-a\left( \varphi +\pi ,\theta \right)
\hspace{5mm}&\text{ for all }\left( \varphi ,\theta \right) \in \mathbb{R}^2,
\label{a1} \\
\int_{0}^{\pi }a\left( s,\theta \right) \binom{\cos s}{\sin s}\ ds=\binom{0}{%
0} \ &\text{ for all }\theta \in \mathbb{R},  \label{a2}
\end{align}
Let $\boldsymbol{V}$ and $\boldsymbol{W}$ be as in (\ref{movbas}) and
suppose that $h:(0,\pi)\times\mathbb{R}\to\mathbb{R}$ is defined by:
\begin{equation}
h\left( \varphi ,\theta \right) :=-\frac{\int_{0}^{\varphi }\sin \left(
\varphi -s\right) ~\partial _{\theta }a(s,\theta )~ds}{\sin \varphi }.
\label{ha}
\end{equation}

\begin{enumerate}
\item Then the definition in (\ref{ha}) can be continuously extended to $%
\mathbb{R}^{2}$. The extended $h$ is such that
\begin{align}
h(\varphi ,\theta )&=h(\varphi+\pi ,\theta )\ \text{ for all }(\varphi ,\theta )
\in \mathbb{R}^2,\label{gdh}\\
h(\varphi ,\theta )&=0\ \text{ for all }(\varphi ,\theta )\in \pi \mathbb{Z}%
\times \mathbb{R} \label{hala}
\end{align}
and satisfies
\begin{equation}
(\varphi ,\theta )\mapsto h(\varphi ,\theta )\boldsymbol{W}(\theta )\in C_{%
\mathrm{p,c}}^{1}(\mathbb{R}^{2}).  \label{holo}
\end{equation}%

\item There exists
\begin{equation}
r_{0}(a)\in \left[ \left\Vert a\right\Vert _{\infty },\left\Vert
a\right\Vert _{\infty }+\left\Vert \partial _{\theta }a\right\Vert _{\infty
}+\left\Vert \partial _{\theta }^{2}a\right\Vert _{\infty }\rule{0cm}{0.4cm}%
\right] ,  \label{rah}
\end{equation}%
such that for all $r\geq r_{0}(a)$ and $\boldsymbol{X}_{0}\in \mathbb{R}^{3}$%
, the surface $\boldsymbol{X}(\mathsf{S})$, defined by
\begin{equation}
\boldsymbol{X}(\varphi ,\theta )=\boldsymbol{X}_{0}+\int_{0}^{\varphi
}\left( r-a\left( s,\theta \right) \right) \boldsymbol{V}(s,\theta
)\,ds+h\left( \varphi ,\theta \right) \boldsymbol{W}(\theta ),
\label{3dform}
\end{equation}%
describes the boundary of a body of constant width $2r$.

\item Moreover, with $a(\cdot ,\cdot )$ as above, the function $h$ in (\ref%
{ha}) is the unique possibility in order that $\boldsymbol{X}$ in (\ref%
{3dform}) describes the boundary of a body of constant width $2r$.
\end{enumerate}
\end{theorem}

\begin{remark}
Although $%
a\in C_{\mathrm{p}}^{2}(\mathbb{R}^2)$ will imply that (\ref{holo}) holds,
one finds at most $\boldsymbol{X}\in C^{0,1}_{\mathrm{p,c}}(\mathbb{R}^2)$.
Hence the induced parametrization $\mathbb{S}^{2}\rightarrow \partial G$ is not necessarily
a diffeomorphism. It will only be the  a diffeomorphism for $r>r_{0}(a)$ and in general
not for $r=r_{0}(a)$. By taking $r>r_0(a)$ one obtains a $C^{1,1}$-surface
with a distance $\varepsilon=r-r_0(a)$ from the body of constant width for $%
r=-r_0(a)$ where Lipschitz is optimal. The surface for $r>r_0(a)$ will also be
a boundary for a body of constant width. Our construction will be illustrated by an example
in Section 4. There the example has $a\in C_{\mathrm{p}}^{1,1}(\mathbb{R}^2)$ and is such
that $r=1.25348\approx r_{0}(a)$, and for  $r=r_{0}(a)$ the surface $\boldsymbol{X}$ will
not be a diffeomorphism. The value of $r_0(a)$ can be computed numerically by finding the smallest
$r\ge\left\|a\right\|_\infty$ such that $T(r_0)\ge 0$ in (\ref{T}) and  $D(r_0,\varphi,\theta)\ge 0$ in (\ref{Det})
holds for all $\varphi,\theta$. Notice both are parabola in $r$ with minima before $\left\|a\right\|_\infty$.
\end{remark}

We have assumed that $a\in C_{\mathrm{p}}^{2}(\mathbb{R}^{2})$, which is
sufficient for describing a 3d set of constant width for $r$ large, but
certainly more than necessary for $h$ and $\boldsymbol{X}$ to be
well-defined. Necessary for $h$ to be well-defined will be $L^{\infty }$
bounds for $a,a_{\theta }$ and $a_{\theta \theta }$. For the 2d case a
necessary and sufficient restriction appears, namely $r\geq
r_{0}(a):=\left\Vert a\right\Vert _{\infty }$. In 3d this condition is still
necessary but not sufficient. To have a differentiable parametrization in 3d
a bound appears that contains $\partial _{\theta }h$. We are however not
able to quantify such a bound more precisely like in 2d.\medskip

We first prove some results for $h$ that we gather in the next lemma. In fact,
Lemma \ref{hlem} contains the first item of Theorem \ref{recipe3d}.

\begin{lemma}
\label{hlem}Suppose that $a\in C_{\mathrm{p}}^{2}(\mathbb{R}^{2})$ satisfies
(\ref{a1}) and (\ref{a2}). Then $h$ in (\ref{ha}) can be continuously extended to $%
\mathbb{R}^{2}$ such that (\ref{ha}) holds for all $\varphi\in\mathbb{R}\setminus\pi\mathbb{Z}$
and $\theta\in\mathbb{R}$.
Moreover, one finds that (\ref{holo}) holds, that equality (\ref{gdh}) holds:
\begin{equation*}
h\left( \varphi ,\theta \right) = h\left( \varphi +\pi ,\theta \right)
 \text{ for all }(\varphi ,\theta)\in \mathbb{R}^{2}
\end{equation*}
and the following estimates: for all $(\varphi
,\theta )\in \mathbb{R}^{2}$
\begin{align}
\left\vert h\left( \varphi ,\theta \right) \right\vert &\leq \left\Vert
a_{\theta }\right\Vert _{\infty }\left\vert \sin \varphi \right\vert ,
\label{tricky} \\
\left\vert h_{\varphi }\left( \varphi ,\theta \right) \right\vert &\leq
\left\Vert a_{\theta }\right\Vert _{\infty },  \label{trickyphi} \\
\left\vert h_{\theta }\left( \varphi ,\theta \right) \right\vert &\leq
\left\Vert a_{\theta \theta }\right\Vert _{\infty }\left\vert \sin \varphi
\right\vert .  \label{trickytheta}
\end{align}
\end{lemma}

\begin{proof}
A priori $h$ is defined for $\varphi \in \left( 0,\pi \right) $ with the periodicity
in the $\theta$-direction being a consequence of the assumption that $a$ satisfies
(\ref{perio1}) and (\ref{perio2}). To consider the extension in the $\varphi$-direction
first let us focus on the enumerator
for $h$ in formula (\ref{ha}). The enumerator is $C^1$, since $a$ is $C^2$. By (\ref%
{a2}) we find that
\begin{equation}
\int_{0}^{\pi }\sin \left( \varphi -s\right) ~a_{\theta }(s,\theta
)~ds=\partial _{\theta }\left( \int_{0}^{\pi }\sin \left( \varphi -s\right)
~a(s,\theta )~ds\right) =0.  \label{hah}
\end{equation}
Moreover, for $\varphi \in \left[ 0,\frac{1}{2}\pi \right] $ we use $1\leq 1+\cos\varphi$ and find%
\begin{multline*}
\left\vert \int_{0}^{\varphi }\sin \left( \varphi -s\right) ~\partial
_{\theta }a(s,\theta )~ds\right\vert \leq \left\Vert \partial _{\theta
}a\right\Vert _{\infty }\int_{0}^{\varphi }\sin \left( \varphi -s\right) ds
\\
=\left\Vert \partial _{\theta }a\right\Vert _{\infty }\left( 1-\cos \varphi
\right) \leq \left\Vert \partial _{\theta }a\right\Vert _{\infty }\left(
1-\cos \varphi \right) \left( 1+\cos \varphi \right) =\left\Vert \partial
_{\theta }a\right\Vert _{\infty }\left\vert \sin \varphi \right\vert ^{2},
\end{multline*}%
while for $\varphi \in \left( \frac{1}{2}\pi ,\pi \right] $ we use $1\leq 1-\cos\varphi$ and (\ref{hah}) to obtain that
\begin{multline*}
\left\vert \int_{0}^{\varphi }\sin \left( \varphi -s\right) ~\partial
_{\theta }a(s,\theta )~ds\right\vert = \left\vert \int_{\varphi}^{\pi }\sin \left( \varphi -s\right) ~\partial
_{\theta }a(s,\theta )~ds\right\vert  \leq \left\Vert \partial _{\theta
}a\right\Vert _{\infty }\left\vert \int_{\varphi }^{\pi }\sin \left( \varphi
-s\right) ds\right\vert \\
=\left\Vert \partial _{\theta }a\right\Vert _{\infty }\left( 1+\cos \varphi
\right) \leq \left\Vert \partial _{\theta }a\right\Vert _{\infty }\left(
1-\cos \varphi \right) \left( 1+\cos \varphi \right) =\left\Vert \partial
_{\theta }a\right\Vert _{\infty }\left\vert \sin \varphi \right\vert ^{2}.
\end{multline*}%
Hence we find that $h$ can be continuously extended by $0$ for $\varphi\in\left\{0,\pi\right\}$ and
\begin{equation*}
\left\vert h\left( \varphi ,\theta \right) \right\vert \leq \left\Vert
\partial _{\theta }a\right\Vert _{\infty }\left\vert \sin \varphi
\right\vert ,
\end{equation*}%
which is (\ref{tricky}) at least on $\left[ 0,\pi \right] \times \mathbb{R}$.

Taking the formula in (\ref{ha}) for $\varphi+\pi\in(\pi,2\pi)$, we find by
(\ref{a2}), a substitution and (\ref{a1})
\begin{align*}
h\left( \varphi +\pi ,\theta \right) & =-\frac{\int_{0}^{\varphi +\pi }\sin
\left( \varphi +\pi -s\right) \partial _{\theta }a\left( s,\theta \right) ds%
}{\sin \left( \varphi +\pi \right) }=\frac{\int_{\pi }^{\varphi +\pi }\sin
\left( \varphi +\pi -s\right) \partial _{\theta }a\left( s,\theta \right) ds%
}{\sin \varphi } \\
& =\frac{\int_{0}^{\varphi }\sin \left( \varphi -s\right) \partial _{\theta
}a\left( s-\pi ,\theta \right) ds}{\sin \varphi }=-\frac{\int_{0}^{\varphi
}\sin \left( \varphi -s\right) \partial _{\theta }a\left( s,\theta \right) ds%
}{\sin \varphi }=h\left( \varphi ,\theta \right)
\end{align*}%
and we find that the definition of $h$ is well-defined on $\left( \pi,2\pi \right)
\times \mathbb{R}$ and at least there (\ref{gdh}) holds. Then $h$ can be
extended continuously by $0$ and (\ref{tricky}) holds
for $\varphi \in \left\{ 0,\pi ,2\pi \right\} $. This allows us to use the
definition in (\ref{ha}) for $h$ for all $\varphi $ with $\sin \varphi \neq
0 $ and to set $h=0$ whenever $\sin \varphi =0$. Moreover, since $a_{\theta
}\in C_{\mathrm{p}}^{1}(\mathbb{R}^{2})$ we find that $h$ and also $h\,%
\boldsymbol{W}$ satisfies (\ref{perio1}). One also finds that (\ref{hala})
holds true.

For (\ref{perio2}) note that%
\begin{align*}
h(-\varphi ,\theta +\pi )& =-\frac{\int_{0}^{-\varphi }\sin \left( -\varphi
-s\right) ~\partial _{\theta }a(s,\theta +\pi )~ds}{\sin \left( -\varphi
\right) }=\frac{\int_{0}^{-\varphi }\sin \left( -\varphi -s\right) ~\partial
_{\theta }a(-s,\theta )~ds}{\sin \varphi } \\
& =-\frac{\int_{0}^{\varphi }\sin \left( -\varphi +s\right) ~\partial
_{\theta }a(s,\theta )~ds}{\sin \varphi }=\frac{\int_{0}^{\varphi }\sin
\left( \varphi -s\right) ~\partial _{\theta }a(s,\theta )~ds}{\sin \varphi }%
=-h(\varphi ,\theta )
\end{align*}%
and with%
\begin{equation*}
h(-\varphi ,\theta +\pi )\boldsymbol{W}(\theta +\pi )=-h(\varphi ,\theta )%
\boldsymbol{W}(\theta +\pi )=h(\varphi ,\theta )\boldsymbol{W}(\theta )
\end{equation*}%
we indeed find (\ref{perio2}) for $h\boldsymbol{W}$. Moreover, $a_{\theta }\in
C_{\mathrm{p}}^{1}(\mathbb{R}^{2})$ implies $h\in C^{1}(%
\mathbb{R}^{2})$ and hence (\ref{holo}) holds.

Since the $\theta $-dependence only comes through $a$ the estimate in (\ref%
{trickytheta}) is proven similarly as for (\ref{tricky}). For (%
\ref{trickyphi}) we use a straightforward computation from (\ref{ha}) and using (\ref{a2}) to find
\begin{equation}
h_{\varphi }\left( \varphi ,\theta \right) =-\frac{\int_{0}^{\varphi }\sin
s~\partial _{\theta }a(s,\theta )~ds}{\left( \sin \varphi \right) ^{2}}=\frac{\int_{\varphi}^{\pi }\sin
s~\partial _{\theta }a(s,\theta )~ds}{\left( \sin \varphi \right) ^{2}}.\label{hafie}
\end{equation}%
It is sufficient to prove the estimates for $\varphi\in[0,\pi]$. So we proceed
for $\varphi\in[0,\frac12\pi]$  by
\begin{equation*}
  \left\vert h_{\varphi }\left( \varphi ,\theta \right) \right\vert \leq
\frac{\int_{0}^{\varphi }\sin s~ds\left\Vert a_{\theta }\right\Vert _{\infty }
 }{\left( \sin \varphi \right) ^{2}}=\frac{1-\cos\varphi}{(\sin \varphi)^2}\left\Vert a_{\theta
 }\right\Vert _{\infty }\leq\frac{(1-\cos\varphi)(1+\cos\varphi)}{(\sin \varphi)^2}\left\Vert a_{\theta
 }\right\Vert _{\infty }=\left\Vert a_{\theta }\right\Vert _{\infty }
\end{equation*}
and for  $\varphi\in[\frac12\pi,\pi]$  by
\begin{equation*}
  \left\vert h_{\varphi }\left( \varphi ,\theta \right) \right\vert \leq
\frac{\int_{\varphi}^{\pi }\sin s~ds\left\Vert a_{\theta }\right\Vert _{\infty }
 }{\left( \sin \varphi \right) ^{2}}=\frac{1+\cos\varphi}{(\sin \varphi)^2}\left\Vert a_{\theta
 }\right\Vert _{\infty }\leq\frac{(1-\cos\varphi)(1+\cos\varphi)}{(\sin \varphi)^2}\left\Vert a_{\theta
 }\right\Vert _{\infty }=\left\Vert a_{\theta }\right\Vert _{\infty }.
\end{equation*}
The estimate in (\ref{trickytheta}) follows as the one in (\ref{tricky}), which concludes the proof of Lemma \ref{hlem}.
\end{proof}

Proofs of Theorem \ref{recipe3d} and of the converse result in the next theorem
are given in Section \ref{sec5}.\pagebreak

\begin{theorem}[All bodies of constant width are represented by (\protect\ref%
{3dform})]
\label{all}\mbox{}
\begin{enumerate}
\item Each body of constant width is described by (\ref{3dform}) for some $%
a\in L_\mathrm{p}^{\infty }(\mathbb{R}^2)$ with $\theta \mapsto
a(\varphi,\theta)$ uniformly Lipschitz on $\mathbb{R}\times\mathbb{R}$ and $a(\cdot,\cdot)$
satisfying (\ref{a1}) and (\ref{a2}), with some $r\geq \left\Vert
a\right\Vert _{L^{\infty }(\mathsf{S})}$ and with $h$ given by
\begin{equation}
h\left( \varphi ,\theta \right) =\frac{\displaystyle\lim_{\varepsilon \rightarrow
0}\int_{0}^{\varphi }\frac{a(s,\theta )-a(s,\theta +\varepsilon )}{%
\varepsilon }\sin \left( \varphi -s\right) ds}{\sin \varphi }.
\label{first-h}
\end{equation}

\item Concerning regularity we have
\begin{equation}
\left( \varphi ,\theta \right)\mapsto h\left( \varphi ,\theta \right) \left(
\begin{array}{c}
-\sin \theta \\
\cos \theta%
\end{array}%
\right) \in C^{0,1}_{\mathrm{p,c}}\left( \mathbb{R}^2\right)\label{hLip}
\end{equation}
and $h$ satisfies (\ref{gdh}) and (\ref{hala}):
\begin{align*}
&h\left( \varphi+\pi ,\theta \right)=h\left( \varphi ,\theta \right) \text{ for all }
( \varphi ,\theta )\in\mathbb{R}\times\mathbb{R},\\
&h\left( \varphi ,\theta \right)=0 \hspace{1cm} \text{ for all }
( \varphi ,\theta )\in\pi\mathbb{Z}\times\mathbb{R}.
\end{align*}%
Moreover, if $a(\cdot,\cdot)$ is such that
\begin{enumerate}
\item[i.] $a,\partial _{\theta }a\in C_\mathrm{p}^{0}(\mathbb{R}^2)$, then
 (\ref{ha}) holds true;

\item[ii.] $a,\partial _{\theta }a\in C_{\mathrm{p}}^{1}(\mathbb{R}^2)$,
then  (\ref{holo}) holds true.
\end{enumerate}
\end{enumerate}
\end{theorem}

\section{An example}

\begin{figure}[H]
\centering
\includegraphics[height=.25\textwidth]{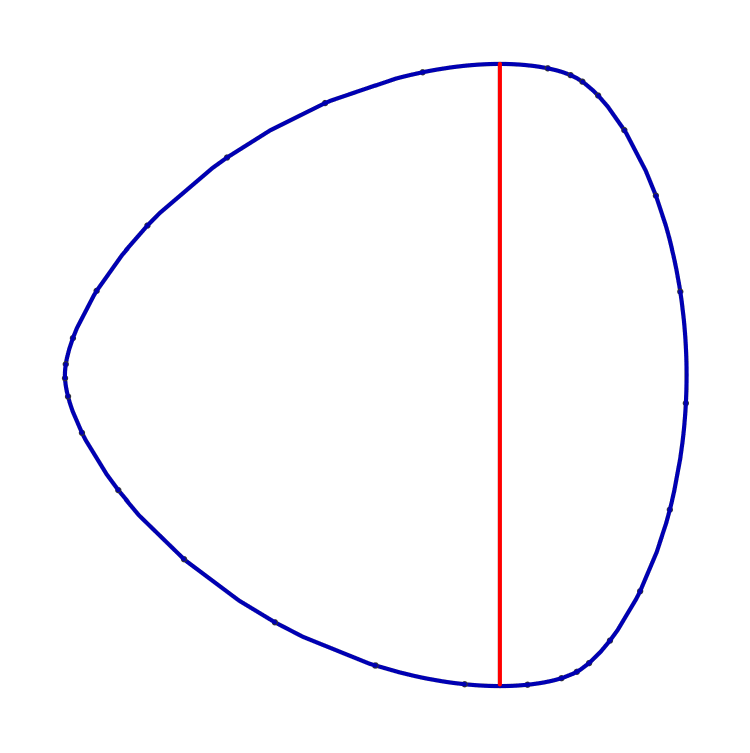}\hspace{6mm} \includegraphics[height=.255%
\textwidth]{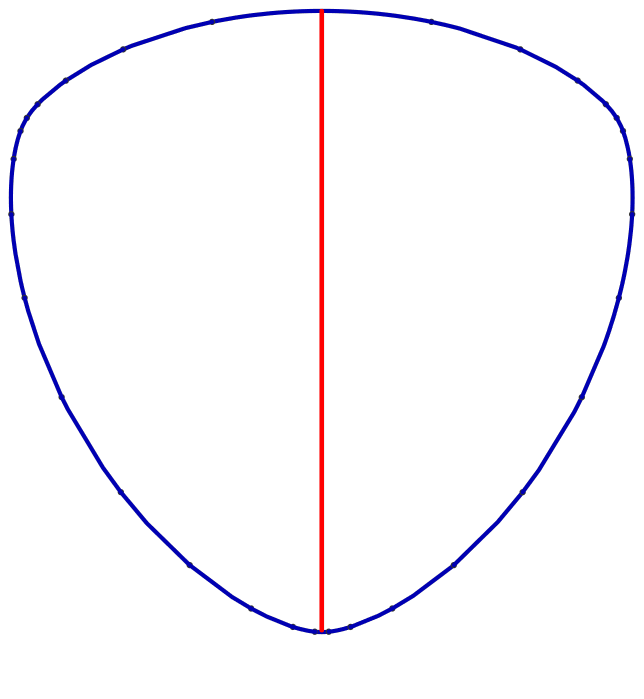}\hspace{2mm} \includegraphics[height=.3%
\textwidth]{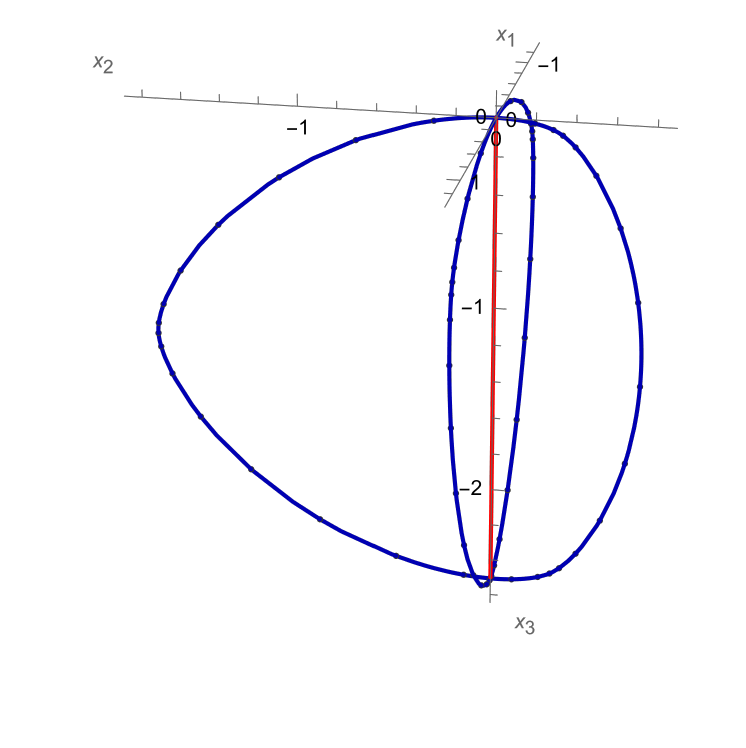}
\caption{The 2d sets with $a_1$ and $a_2$ for $\boldmath{x}(0)$ on top, and the
combination in 3d by perpendicular planes and joining the red axes.}
\label{2dbild}
\end{figure}

The formulas are rather technical and in order to illustrate that
(\ref{3dform}) does deliver a body of constant width, we give an actual
construction in a case that is computable. The example shows a body of
constant width connecting two triangular 2d-domains of constant width based
on the 2d-formula. In addition to $\boldsymbol{x}_{0}=(0,0)$ and $r=1$ we
use in Fig.~\ref{2dbild}:

\begin{itemize}
\item for the figure on the left: $a(s)=a_{1}(s):=-\cos (3s)$;

\item for the figure in the middle: $a(s)=a_{2}(s):=\sin (3s)$.
\end{itemize}
One directly checks that conditions (\ref{a1-2d}) and (\ref{a2-2d})
are satisfied for $a_1$ and $a_2$.\medskip

The object on the right of Fig.~\ref{2dbild} combines these two curves in a
3d-setting in orthogonal planes with the red line as common intersection. In
order to find a smooth perturbation from the horizontal to the vertical
curve by curves whose projections will be 2d-curves of constant width $1$,
we use the following:%
\begin{equation}
a\left( \varphi ,\theta \right) :=\left( \cos \theta \right)
^{2}a_{1}(\varphi )+\left| \sin \theta \right| \sin \theta~a_{2}(\varphi ).
\label{combiaa}
\end{equation}

Since $a_1$ and $a_2$ satisfy (\ref{a1-2d}) and (\ref{a2-2d}) it follows
that $a$ defined in (\ref{combiaa}) satisfies (\ref{a1}) and (\ref{a2}). The
periodicity condition in (\ref{perio1}) one checks directly; the one in (\ref{perio2}) follows
since $a_1(\varphi)=a_1(-\varphi)$ and $a_2(\varphi)=-a_2(-\varphi)$.
On may compute $\left\|a\right\|_\infty= 1$, $\left\|\partial_\theta a\right\|_\infty=\sqrt{2}$
and $\left\|\partial_\theta^2 a\right\|_\infty= 2\sqrt{2}$.

One has $r_0(a)\geq \max\left(r_0(a_1),r_0(a_2)\right)=1$. Concerning the value of
$r_0(a)$ for $a$ in (\ref{combiaa}) a numerical estimate for the expression in (\ref{Det}) to be positive
shows $r_0(a)\approx 1.25348$, which lies inside the interval given in (\ref{rah}).

\begin{figure}[H]
\centering
\includegraphics[width=.42\textwidth]{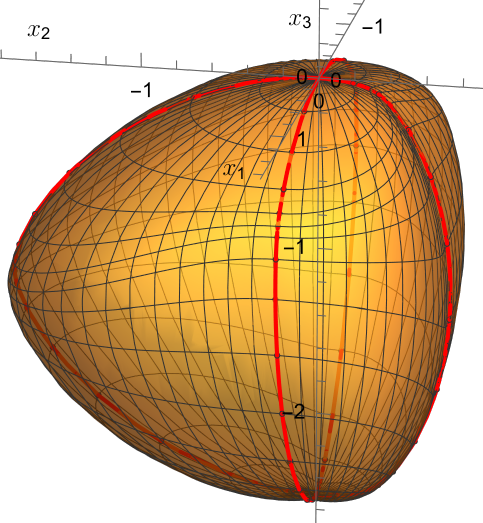}\hspace{3mm}\includegraphics[width=.42\textwidth]{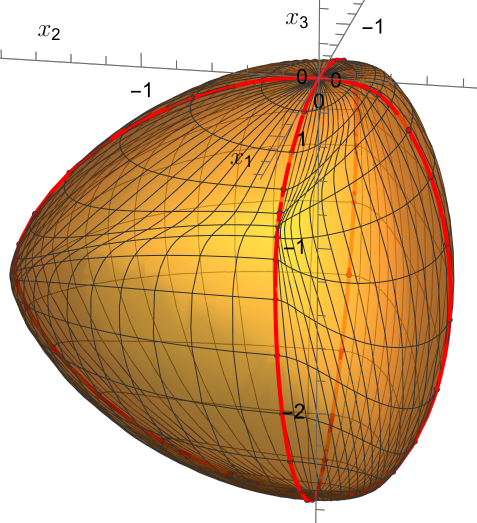}
\caption{On the left the intermediate construction still without the $h$.
It consists of a rotating family of
2d-sets of constant width for each $\protect\theta$. This is not a body of
constant width and not even convex. On the right is the corresponding body
of constant width as the final result
with the shift by $h$ in the direction $\boldsymbol{\Psi}$ from (\protect\ref%
{basis}). Both red curves originate from the curves from Fig.~\protect\ref{2dbild}.
The surface on the right does not look smooth everywhere
 as indeed here $r=1.25348$, the numerical approximation of $r_0(a)$. These red curves give two planar
 curves of constant width, since $\partial_\theta a(\varphi,\theta)=0$ for
 all $\left(\varphi,\theta\right)\in\mathbb{R}\times\pi\mathbb{Z}$.}
\label{3dbild}
\end{figure}

The $a$ in (\ref{combiaa}) is used to produce the sketch on the left in Fig.~%
\ref{3dbild} using the formula in (\ref{3dform}) without the $h$-term. Each
intersection with a plane containing the vertical (red) line $%
\left\{\lambda (0,0,1); \lambda\in\mathbb{R}\right\}$ will produce a 2d set
of constant width. After the modification with the additional $h$-term in
(\ref{3dform}) does one indeed find a 3d set of constant width, which is found
on the right of Fig.~\ref{3dbild}.\medskip

Although each body of constant width can be constructed through the formula
in \ref{3dform} it is relatively easy if one connect two curves of constant
width as  \ref{combiaa}. More examples can be found in Fig.~\ref{exes}.

\section{Proofs of the two theorems}\label{sec5}

For the standard inner product of $u,v\in \mathbb{R}^{n}$ we use $%
\left\langle u,v\right\rangle $. The notation $u\cdot v$ is used for
componentwise multiplication, which includes but can be more general than
the inner product. Let us start by introducing three vectors for a more
concise notation:
\begin{equation}
\boldsymbol{\Theta }=\left(
\begin{array}{c}
\cos \theta \\
\sin \theta \\
0%
\end{array}%
\right) \text{, }\boldsymbol{\Psi }=\left(
\begin{array}{c}
-\sin \theta \\
\cos \theta \\
0%
\end{array}%
\right) \text{ \ and }\boldsymbol{\Xi }=\left(
\begin{array}{c}
0 \\
0 \\
1%
\end{array}%
\right) .  \label{basis}
\end{equation}%
These three directions constitute a $\theta $-dependent orthonormal basis in
$\mathbb{R}^{3}$ that turns out to be convenient for our parametrization.
The following identities hold true:%
\begin{equation}
\partial _{\theta }\boldsymbol{\Theta }=\boldsymbol{\Psi }\text{ and }%
\partial _{\theta }\boldsymbol{\Psi }=-\boldsymbol{\Theta .}  \label{iden}
\end{equation}%
Also note that our initial basis (\ref{movbas}) can be expressed in term of (\ref{basis}):%
\begin{align}
\boldsymbol{U}\left( \varphi ,\theta \right) & =\left(
\begin{array}{c}
\sin \varphi \cos \theta \\
\sin \varphi \sin \theta \\
\cos \varphi%
\end{array}%
\right) =\cos \varphi ~\boldsymbol{\Xi }+\sin \varphi ~\boldsymbol{\Theta }%
=\left(
\begin{array}{c}
\sin \varphi \\
\cos \varphi%
\end{array}%
\right) \cdot \left(
\begin{array}{c}
\boldsymbol{\Theta } \\
\boldsymbol{\Xi }%
\end{array}%
\right) ,  \label{notp} \\
\boldsymbol{V}\left( \varphi ,\theta \right) & =\left(
\begin{array}{c}
\cos \varphi \\
-\sin \varphi%
\end{array}%
\right) \cdot \left(
\begin{array}{c}
\boldsymbol{\Theta } \\
\boldsymbol{\Xi }%
\end{array}%
\right) \text{ \ and \ }\boldsymbol{W}\left( \theta \right) =\boldsymbol{%
\Psi }.  \label{notp2}
\end{align}
The dot product $\cdot $ in (\ref{notp}), (\ref{notp2}) is a more convenient  notation in the following
proofs.\medskip

\begin{proof}[Proof of Theorem \protect\ref{recipe3d}]
We will have to show that $\boldsymbol{X}$ in (\ref{3dform}) is a regular
parametrization and secondly, that the resulting surface will yield a body
of constant width. For both aspects we need to consider $\partial _{\varphi }%
\boldsymbol{X}\left( \varphi ,\theta \right) $ and $\partial _{\theta }%
\boldsymbol{X}\left( \varphi ,\theta \right) $.\medskip

$\blacktriangleright $ \emph{Computation of $\partial _{\varphi }\boldsymbol{%
X}\left( \varphi ,\theta \right)$ and $\partial _{\theta }\boldsymbol{X}\left( \varphi ,\theta \right)$.} We will check first that $%
\boldsymbol{X}$ in (\ref{3dform}) is a regular parametrization of the
boundary $\partial G$ of a body of constant width for $r$ large enough, that
is
\begin{equation}
\boldsymbol{\tilde{X}}:\mathbb{S}^{2}\rightarrow \partial G\text{ defined by
}\boldsymbol{\tilde{X}}\left( \omega \right) :=\boldsymbol{X}\left( \varphi
,\theta \right) \text{ for }\omega =\boldsymbol{U}(\varphi ,\theta )
\label{XtildeX}
\end{equation}%
is $C^{1}$, one-to-one and onto, and even a diffeomorphism. With the
notation from (\ref{basis}) we can rewrite (\ref{3dform}) as%
\begin{equation}
\boldsymbol{X}\left( \varphi ,\theta \right) =\boldsymbol{X}%
_{0}+\int_{0}^{\varphi }\left( r-a(s,\theta )\right) \left(
\begin{array}{c}
-\sin s \\
\cos s%
\end{array}%
\right) ds\cdot \left(
\begin{array}{c}
\boldsymbol{\Xi } \\
\boldsymbol{\Theta }%
\end{array}%
\right) +h(\varphi ,\theta )\ \boldsymbol{\Psi .}  \label{Xshort}
\end{equation}%
One computes that
\begin{equation}
\partial _{\varphi }\boldsymbol{X}\left( \varphi ,\theta \right) =\left(
r-a(\varphi ,\theta )\right) \left(
\begin{array}{c}
-\sin \varphi \\
\cos \varphi%
\end{array}%
\right) \cdot \left(
\begin{array}{c}
\boldsymbol{\Xi } \\
\boldsymbol{\Theta }%
\end{array}%
\right) +\partial _{\varphi }h(\varphi ,\theta )\ \boldsymbol{\Psi }\label{split-Xphi}
\end{equation}%
and that
\begin{multline}
\partial _{\theta }\boldsymbol{X}\left( \varphi ,\theta \right)
=-\int_{0}^{\varphi }\partial _{\theta }a(s,\theta )\left(
\begin{array}{c}
-\sin s \\
\cos s%
\end{array}%
\right) ds\cdot \left(
\begin{array}{c}
\boldsymbol{\Xi } \\
\boldsymbol{\Theta }%
\end{array}%
\right) -h(\varphi ,\theta )\ \boldsymbol{\Theta }~+ \\
\left( \int_{0}^{\varphi }\left( r-a(s,\theta )\right) \cos s~ds+\partial
_{\theta }h(\varphi ,\theta )\right) \boldsymbol{\Psi }.\label{split-Xthe}\medskip
\end{multline}

$\blacktriangleright $ \emph{Invariant normal direction.} The next step is
to show that the outward normal direction at $\boldsymbol{X}%
\left( \varphi ,\theta \right) $ satisfies:
\begin{equation}
\nu_{\boldsymbol{\tilde{X}}\left( \omega \right)}=\pm\omega \text{ for all }%
\omega \in \mathbb{S}^{2}  \label{pregm}
\end{equation}%
where $\omega=\boldsymbol{U}(\varphi,\theta)$ as in (\ref{defpar}). Indeed,
we will first show that $\omega$ is perpendicular to $\partial _{\varphi }\boldsymbol{X}\left( \varphi ,\theta
\right) $ and $\partial _{\theta }\boldsymbol{X}\left( \varphi ,\theta \right)$.
Taking $h$ as in (\ref{ha})
is in fact the only possible choice such that
\begin{equation}
\omega \cdot \partial _{\theta }\boldsymbol{X}\left( \varphi ,\theta \right) =0  \label{derig}
\end{equation}
holds. Indeed with this $h$ we may rewrite (\ref{split-Xthe}), at least when $\sin\varphi\ne 0$, to
\begin{multline}
\partial _{\theta }\boldsymbol{X}\left( \varphi ,\theta \right)
=\int_{0}^{\varphi }\partial _{\theta }a(s,\theta )\sin s~ds~\boldsymbol{\Xi
}-\frac{\cos \varphi }{\sin \varphi }\int_{0}^{\varphi }\partial _{\theta
}a(s,\theta )
\sin s~ds\ \boldsymbol{\Theta}~+ \\
\left( \int_{0}^{\varphi }\left( r-a(s,\theta )\right) \cos s~ds+\partial
_{\theta }h(\varphi ,\theta )\right) \boldsymbol{\Psi}\label{split-Xthe-bis}
\end{multline}
and using $\omega =\cos \varphi \ \boldsymbol{\Xi }+\sin \varphi
\boldsymbol{\ \Theta }$ one directly finds (\ref{derig}).
As in the proof of Lemma \ref{hlem} the factor in front of $\boldsymbol{\Xi}$ in
(\ref{split-Xthe-bis}) can be continuously extended by $0$ when $\sin\varphi=0$.
From (\ref{split-Xphi}) and again with $\omega =\cos \varphi \ \boldsymbol{\Xi }+\sin \varphi
\boldsymbol{\ \Theta }$ we find%
\begin{equation*}
\omega \cdot \partial _{\varphi }\boldsymbol{X}\left( \varphi ,\theta
\right) =0.  \label{werig}
\end{equation*}
If $\partial _{\theta }\boldsymbol{X}$ or $\partial _{\varphi }\boldsymbol{X}$ is trivial
for some $(r,\theta,\varphi)$ one may consider $(r+\varepsilon,\theta,\varphi)$ and
find from (\ref{split-Xphi}) and (\ref{split-Xthe-bis}) that for
$\varepsilon>0$ the corresponding expressions will be nontrivial and (\ref{derig})
will hold for that $(\varphi,\theta)$.

After the homotopy to the sphere furtheron, one may conclude that $\omega$ is the outside normal
for $(r+\varepsilon,\theta,\varphi)$  with $\varepsilon$ large and, by continuity, is an outside normal for
$\boldsymbol{X}(\mathbb{R}^2)$ at $(\theta,\varphi)$ when $r> r_0(a)$ with $r_0(a)$ to be defined in (\ref{rzero}).\medskip

$\blacktriangleright $ \emph{Well defined parametrization.} For $a\in C_%
\mathrm{p}^{2}( \mathbb{R}^2) $ Lemma \ref{hlem} implies that $h$ is
well-defined and $h\boldsymbol{W}$ lies in $C_\mathrm{p,c}^{1}(\mathbb{R}^2)$%
. So with (\ref{a1}) and (\ref{a2}) also the expression in (\ref{Xshort})
lies in $C^{1}_\mathrm{p,c}(\mathbb{R}^2)$.

In order to have a regular parametrization it is sufficient that:

\begin{itemize}
\item $\partial _{\varphi }\boldsymbol{X}\times \partial _{\theta }%
\boldsymbol{X}$ is nontrivial on $\left\{ \left( \varphi ,\theta \right) \in
\mathbb{R}^{2};\varphi \not\in \pi \mathbb{Z}\right\} $, and

\item $\partial _{\varphi }\boldsymbol{X}\left( \varphi ,0\right) \times
\partial _{\varphi }\boldsymbol{X}\left( \varphi ,\frac{1}{2}\pi \right) $
is nontrivial for $\varphi \in \left\{ 0,\pi \right\} $.
\end{itemize}

Let us start with the second case for $\varphi =0$, with $\varphi =\pi $
similarly:
\begin{equation*}
\partial _{\varphi }\boldsymbol{X}(\varphi ,0)\times \partial _{\varphi }%
\boldsymbol{X}(\varphi ,\tfrac{1}{2}\pi )=\left(
\begin{array}{c}
r-a(0,0) \\
\partial _{\varphi }h\left( 0,0\right) \\
0%
\end{array}%
\right) \times \left(
\begin{array}{c}
-\partial _{\varphi }h\left( 0,\frac{1}{2}\pi \right) \\
r-a(0,\frac{1}{2}\pi ) \\
0%
\end{array}%
\right) =:\left(
\begin{array}{c}
0 \\
0 \\
T(r)%
\end{array}%
\right) ,
\end{equation*}%
where
\begin{equation}
T(r)=\left( r-a(0,0)\right) \left( r-a(0,\tfrac{1}{2}\pi )\right) +\partial
_{\varphi }h\left( 0,0\right) \partial _{\varphi }h\left( 0,\tfrac{1}{2}\pi
\right) .  \label{T}
\end{equation}%
Since
$\left\vert h_{\varphi }\left( \varphi ,\theta \right) \right\vert \leq
\left\Vert \partial_{\theta }a \right\Vert _{\infty }$ holds, see (\ref{trickyphi}), a sufficient condition for $T(r)>0$ is
\begin{equation*}
r>\left\Vert a\right\Vert _{\infty }+\left\Vert \partial _{\theta
}a\right\Vert _{\infty }.
\end{equation*}
Note that $T^\prime(r)\ge 0$ for $r\ge \left\|a\right\|_\infty$.\medskip

For $\varphi \not\in \pi \mathbb{Z}$, using (\ref{werig}) and (\ref{derig}),
which state that $\omega $ is perpendicular to $\partial _{\varphi }%
\boldsymbol{X}$ and $\partial _{\theta }\boldsymbol{X}$, a simple way of
checking that $\partial _{\varphi }\boldsymbol{X}\times \partial _{\theta }%
\boldsymbol{X}$ is nontrivial, is to show that $\boldsymbol{U}(\varphi,\theta) \cdot \left( \partial
_{\varphi }\boldsymbol{X}\times \partial _{\theta }\boldsymbol{X}\right)
\neq 0$. Using the orthonormal basis $\left\{ \boldsymbol{\Theta },%
\boldsymbol{\Psi ,\Xi }\right\} $ we obtain from (\ref{split-Xphi}) and (\ref{split-Xthe-bis}):
\begin{align}
&\hspace{1cm}D(r,\varphi,\theta):=\boldsymbol{U}(\varphi,\theta) \cdot \left( \partial _{\varphi }\boldsymbol{X}(\varphi
,\theta )\times \partial _{\theta }\boldsymbol{X}(\varphi ,\theta )\rule%
{0mm}{3mm}\right) \notag\\[2mm]
=& \det \left(
\begin{array}{ccc}
\sin \varphi & \left( r-a(\varphi ,\theta )\right) \cos \varphi & -\frac{%
\cos \varphi }{\sin \varphi }\int_{0}^{\varphi }\partial _{\theta
}a(s,\theta )\sin s~ds \\
0 & \partial _{\varphi }h(\varphi ,\theta ) & \int_{0}^{\varphi }\left(
r-a(s,\theta )\right) \cos s~ds+\partial _{\theta }h(\varphi ,\theta ) \\
\cos \varphi & -\left( r-a(\varphi ,\theta )\right) \sin \varphi &
\int_{0}^{\varphi }\partial _{\theta }a(s,\theta )\sin s~ds%
\end{array}%
\right) \notag\\[2mm]
&\hspace{-11mm}=\left( r-a(\varphi ,\theta )\right) \left( \partial _{\theta
}h(\varphi ,\theta )+\int_{0}^{\varphi }\left( r-a(s,\theta )\right) \cos
s~ds\right) +\partial _{\varphi }h(\varphi ,\theta )\frac{\int_{0}^{\varphi
}\sin s\ \partial _{\theta }a(s,\theta )~ds}{\sin \varphi } \notag\\
=&\left( r-a(\varphi ,\theta )\right) \left( \int_{0}^{\varphi }\left(
r-a(s,\theta )\right) \cos s~ds+\partial _{\theta }h(\varphi ,\theta
)\right) -\left( \partial _{\varphi }h(\varphi ,\theta )\right) ^{2}\sin
\varphi. \label{Det}
\end{align}
In the last step we used (\ref{hafie}). Assuming $r\geq \left\| a\right\|_\infty$ we have
\begin{align*}
   \text{ for }\varphi\in [0,\tfrac12\pi]:\quad\int_{0}^{\varphi }\left( r-a(s,\theta )\right) \cos s~ds&
   \ge \left(r-\left\|a\right\|_\infty\right)\sin\varphi,\\
   \text{ for }\varphi\in [\tfrac12\pi,\pi]:\quad\int_{0}^{\varphi }\left( r-a(s,\theta )\right) \cos s~ds&
   = -\int_{\varphi }^\pi\left( r-a(s,\theta )\right) \cos s~ds
   \ge \left(r-\left\|a\right\|_\infty\right)\sin\varphi  .
\end{align*}
For $\varphi\in (\pi,2\pi)$ one obtains
similar estimates for $\left|D(r,\varphi,\theta)\right|=-D(r,\varphi,\theta)$. The expression in (\ref{Det})
can now be estimated. Using (\ref{trickyphi}) and
(\ref{trickytheta}) from  Lemma \ref{hlem} we get for $\varphi\not\in\pi\mathbb{Z}$:
\begin{gather*}
\frac{D(r,\varphi,\theta)}{\sin\varphi} \geq \left( \left( r-\left\Vert a\right\Vert _{\infty }\right) ^{2}-\left(
r-\left\Vert a\right\Vert _{\infty }\right) \left\Vert \partial _{\theta
}^{2}a\right\Vert _{\infty }-\left\Vert \partial _{\theta }a\right\Vert
_{\infty }^{2}\right)  \\[2mm]
\geq \left( r-\left\Vert a\right\Vert _{\infty }-\left\Vert \partial
_{\theta }a\right\Vert -\left\Vert \partial _{\theta }^{2}a\right\Vert
_{\infty }\rule{0mm}{4mm}\right) \left( r-\left\Vert a\right\Vert _{\infty
}+\left\Vert \partial _{\theta }a\right\Vert _{\infty }\rule{0mm}{4mm}%
\right)  ,
\end{gather*}%
which is positive whenever
\begin{equation}\label{Deter}
  r\geq \left\Vert a\right\Vert
_{\infty }+\left\Vert \partial _{\theta }a\right\Vert _{\infty }+\left\Vert
\partial _{\theta }^{2}a\right\Vert _{\infty }.
\end{equation}
Moreover, whenever $r\ge \left\|a\right\|_\infty$ (\ref{Det}) also shows
that $\partial_r \left(D(r,\varphi,\theta)/\sin\varphi\right)\ge 0$.

Since $T$ and $\left|D\right|$ for $\varphi\not\in \pi\mathbb{Z}$ are increasing
with respect to $r$ for $r\ge \left\|a\right\|_\infty$, there exists a minimal
\begin{equation}
r_{0}(a)\in \left[ \left\Vert a\right\Vert _{\infty },\left\Vert
a\right\Vert _{\infty }+\left\Vert \partial _{\theta }a\right\Vert _{\infty
}+\left\Vert \partial _{\theta }^{2}a\right\Vert _{\infty }\right]\label{rzero}
\end{equation}%
such that $T$ and $\left|D\right|$ for $\varphi\not\in \pi\mathbb{Z}$ are positive and hence that the parametrization is
well-defined for all $r>r_{0}(a)$. For $r=r_{0}(a)$ the parametrization is
no longer necessarily of class $C^{1}$ or one-to-one. However,
since for all $r>r_{0}(a)$ one will find a body of constant width and all
functions involved are continuous, also the limit by taking $r\downarrow
r_{0}(a)\ne 0$ will give a body of constant width.\medskip

$\blacktriangleright $ \emph{Homotopy to the sphere.} The parametrization is
well-defined for all $r>r_{0}(a)$ and to be able to focus on the dependence
on $r$ we use an explicit $r$ in the following expression (\ref{3dform}) in
this paragraph:
\begin{equation*}
\boldsymbol{X}_{e}\left( r,\varphi ,\theta \right) :=\boldsymbol{X}\left(
\varphi ,\theta \right) \text{ and } \boldsymbol{\tilde{X}}_{e}\left(
r,\omega \right) :=\boldsymbol{\tilde{X}}\left( \omega \right),
\end{equation*}%
with $\boldsymbol{X},\boldsymbol{\tilde{X}}$ from (\ref{3dform}) and (\ref{XtildeX}). We define%
\begin{equation*}
\left( 0,1\right] \times \mathbb{S}^{2}\ni \left( \rho ,\omega \right)
\mapsto \boldsymbol{\tilde{Y}}\left( \rho ,\omega \right) :=%
\rho\,\boldsymbol{\tilde{X}}_{e}\left( \rho ^{-1}r,\omega \right) \in \mathbb{R}%
^{3}.
\end{equation*}%
As one may see from (\ref{Xshort})  one finds that
\begin{equation*}
\boldsymbol{\tilde{Y}}\left( 1,\omega \right) =\boldsymbol{\tilde{%
X}}_{e}\left( r,\omega \right) \text{ and }\boldsymbol{\tilde{Y}}\left(
0,\omega \right) :=\lim_{\rho \downarrow 0}\boldsymbol{\tilde{Y}}\left( \rho
,\omega \right) =r\omega
\end{equation*}%
with all $\boldsymbol{\tilde{Y}}\left( \rho ,\cdot \right) $ for $\rho\in[0,1%
]$ being regular parametrizations and their outward normal directions $\nu$ satisfying as in (\ref{pregm})%
\begin{equation}  \label{normals}
\nu_{\boldsymbol{\tilde{Y}}( \rho ,\omega )} =\omega \text{ for all }\omega
\in \mathbb{S}^{2}.
\end{equation}%
So $\mathbb{R}^{3}\setminus \boldsymbol{\tilde{X}}_{e}\left( r,\mathbb{S}%
^{2}\right) $ has precisely two connected components. We
call $A$ the bounded one.\medskip

$\blacktriangleright $ \emph{Convexity of $A$.} Since the extreme value of $%
\boldsymbol{\tilde{X}}(\mathbb{S}^{2})$ in the direction $\omega $ has
normal $\omega $, and since by (\ref{normals}) $\nu _{\boldsymbol{\tilde{X}}%
\left( \omega \right) }=\omega $, that extreme point is indeed $\boldsymbol{%
\tilde{X}}\left( \omega \right) $. So for each $\omega $ it holds that $%
\boldsymbol{\tilde{X}}(\mathbb{S}^{2})$, except for $\boldsymbol{\tilde{X}}%
\left( \omega \right) $ itself, is on one side of that tangent plane. Hence $%
\bar{A}$ lies on one side of all the tangent planes for $\partial A=%
\boldsymbol{\tilde{X}}(\mathbb{S}^{2})$, which implies that $\bar{A}$ is
convex. See also the proof of Hadamard's Theorem \cite[page 194]{MP}.\medskip

$\blacktriangleright $ \emph{Body of constant width.} According to the
results proved above it is sufficient to show that%
\begin{equation*}
\boldsymbol{\tilde{X}}\left( \omega \right) -\boldsymbol{\tilde{X}}\left(
-\omega \right) =2r\omega \text{ for all }\omega \in \mathbb{S}^{2}.
\end{equation*}%
For the parametrization $\boldsymbol{X}$ with $\boldsymbol{U}$ as in (\ref%
{defpar}) this coincides with%
\begin{equation*}
\boldsymbol{X}(\varphi ,\theta )-\boldsymbol{X}(\varphi +\pi ,\theta )=2r\
\boldsymbol{U}\left( \varphi ,\theta \right) \text{ for all }(\varphi
,\theta )\in \mathsf{S}.
\end{equation*}%
Indeed, using (\ref{Xshort}) we find with (\ref{a1}), (\ref{a2}) and (\ref%
{gdh}) that%
\begin{align*}
& \boldsymbol{X}(\varphi +\pi ,\theta )-\boldsymbol{X}(\varphi ,\theta ) \\
& =\int_{\varphi }^{\varphi +\pi }\left( r-a(s,\theta )\right) \left(
\begin{array}{c}
-\sin s \\
\cos s%
\end{array}%
\right) ds\cdot \left(
\begin{array}{c}
\boldsymbol{\Xi } \\
\boldsymbol{\Theta }%
\end{array}%
\right) +\left( h(\varphi +\pi ,\theta )-h(\varphi ,\theta )\right) \
\boldsymbol{\Psi } \\
& =r\left(
\begin{array}{c}
\cos \left( \varphi +\pi \right) -\cos \varphi \\
\sin \left( \varphi +\pi \right) -\sin \varphi%
\end{array}%
\right) \cdot \left(
\begin{array}{c}
\boldsymbol{\Xi } \\
\boldsymbol{\Theta }%
\end{array}%
\right) =-2r\ \boldsymbol{U}\left( \varphi ,\theta \right) ,
\end{align*}%
as desired.
\end{proof}

\begin{proof}[Proof of Theorem \protect\ref{all}, the derivation of
formula (\ref{3dform}) for some $h$]
Suppose that $G$ is a body of constant width $\boldsymbol{d}_G$. Define $\boldsymbol{X}_{0}\in
\mathbb{R}^{3}$ as the point on $\partial G$ with the largest $x_{3}$%
-coordinate. Since a translation that maps $\boldsymbol{X}_0$ to a fixed point does
not meddle with our arguments, we may assume
\begin{equation}\label{X0}
  \boldsymbol{X}_0=\left(0,0,\tfrac12 \boldsymbol{d}_G\right)^{T}.
\end{equation}

Taking $u=\left( 0,0,1\right) ^{T}$ and $\omega =\left( 0,-\sin
\theta ,\cos \theta \right) ^{T}$ the result of Hadwiger, extended by the
remark of Groemer that bodies of constant width have only regular boundary
points, states that is is sufficient that the projections $P_{\omega }G$\ of
$G$, on each of the planes spanned by $\left\{\boldsymbol{\Theta}(\theta),\boldsymbol{\Xi}\right\}$
with $\theta \in [0,\pi] $, are curves
of constant width $\boldsymbol{d}_{\widehat{P}_{\omega }G}=2r$. Thus by Theorem \ref{Th2}
all those sets can be described by (\ref{2df}) with for each $%
\theta $ some function $a(\cdot,\theta)$ depending on $\theta $ as a
parameter. The value of $r$ is the same for all projections and does not
depend on $\theta $. In other words, a fixed $r$ exists and for each $\theta
$ a mapping $\varphi \mapsto a(\varphi ,\theta )\in L^{\infty }(0,2\pi )$
such that for the corresponding $\boldsymbol{x}$ as in Theorem \ref{recipe2d}
we have, with $\widehat{P}$ as in (\ref{Poms}),
\begin{equation*}
\partial \widehat{P}_{\omega }G=\boldsymbol{x}\left( \left[ 0,2\pi \right]
,\theta \right)
\end{equation*}%
with some $\boldsymbol{x}_0(\theta)=\boldsymbol{x}\left( 0,\theta \right) \in \mathbb{R}^{2}$
in accordance with Theorem \ref{recipe2d} and
\begin{equation*}
\sup \left\{ \left\vert a\left( \varphi ;\theta \right) \right\vert ;0\leq
\varphi \leq \pi \right\} \leq r\ \text{ for all }\theta \in \left[0,\pi
\right] .
\end{equation*} Moreover, the mapping $\varphi \mapsto a(\varphi ,\theta )$
satisfies (\ref{a1-2d}) and (\ref{a2-2d}). Hence (\ref{a1}), (\ref{a2}) and $%
r_{0}(a)\geq \left\Vert a\right\Vert _{L^{\infty }(\mathsf{S})}$ are
necessary conditions.

Since for each $\tilde{\omega}\in \mathbb{S}^{2}$ the set $G$ lies in the
cylinder perpendicular to its projection, in other words, we have $G\subset
P_{\tilde{\omega}}G+\left[ \tilde{\omega}\right] $ with $\left[ \tilde{\omega%
}\right] =\left\{ \lambda \tilde{\omega};\lambda \in \mathbb{R}\right\} $.
It follows that for each%
\begin{equation*}
X_{\ast }\in \partial G\cap \left( \partial P_{\tilde{\omega}}G+\left[
\tilde{\omega}\right] \right)
\end{equation*}%
there is $\left( \varphi ,\theta \right) \in \mathsf{S}$ with $%
\boldsymbol{U}(\varphi ,\theta )\cdot \tilde{\omega}=0$ and a value $h\left(
\varphi ,\theta \right) \in \mathbb{R}$ such that
\begin{equation*}
X_{\ast }=P_{\tilde{\omega}}(X_{\ast })+h(\varphi ,\theta )\,\tilde{\omega},
\end{equation*}%
with $\tilde\omega=\boldsymbol{\Psi}(\theta)$. If $X_{\ast \ast }\in \partial G$
is such that $\left\Vert X_{\ast }-X_{\ast
\ast }\right\Vert =2r$, with $2r$ being the width, also $\left\Vert P_{\tilde{\omega}%
}(X_{\ast })-P_{\tilde{\omega}}(X_{\ast \ast })\right\Vert =2r$ and hence%
\begin{equation*}
X_{\ast \ast }=P_{\tilde{\omega}}(X_{\ast \ast })+h(\varphi ,\theta )\,\tilde{%
\omega},
\end{equation*}%
with the same contribution $h(\varphi ,\theta )\tilde\omega$, which implies that
\begin{equation}
h\left( \varphi ,\theta \right) =h\left( \varphi +\pi ,\theta \right) \text{
for all }\left( \varphi ,\theta \right) \in \mathsf{S}\text{.}  \label{haco}
\end{equation}%
Indeed
\begin{equation}
X_{\ast }=\boldsymbol{X}_0+\boldsymbol{x}(\varphi ,\theta )\cdot \left(
\begin{array}{c}
\boldsymbol{\Xi } \\
\boldsymbol{\Theta}(\theta) %
\end{array}%
\right) +h\left( \varphi ,\theta \right) \, \boldsymbol{\Psi}(\theta).  \label{der}
\end{equation}%
Here (\ref{haco}) follows from the fact that the line through the points of
farthest distance is perpendicular to the plane spanned by
$\left\{\boldsymbol{\Theta}(\theta),\boldsymbol{\Xi}\right\}$. Since for $\varphi \in \left\{ 0,\pi
,2\pi \right\} $ the $X_{\ast }$ in (\ref{der}) does not depend on $\theta $%
, one finds for all $\theta \in [0,\pi]$,
that
\begin{equation*}
h\left( 0,\theta \right) =h\left( \pi ,\theta \right) =h\left( 2\pi ,\theta
\right) =0
\end{equation*}%
and
\begin{equation*}
\boldsymbol{x}(0,0)=\boldsymbol{x}(0,\theta )=\boldsymbol{x}(2\pi ,\theta )=%
\boldsymbol{x}(\pi ,\theta )+\left(
\begin{array}{c}
2r \\
0%
\end{array}%
\right) .
\end{equation*}

The first factor on the right in (\ref{der}) inherits the conditions of the
two-dimensional formula and so for each $\theta \in [0,\pi]$ one finds
$\varphi \mapsto \boldsymbol{x}(\varphi ,\theta )$
as in (\ref{2df}). The formula in (\ref{der})  describes through $%
\boldsymbol{X}:\mathsf{S}\rightarrow \mathbb{R}^{3}$ all points of $\partial G$ by
\begin{equation}
\boldsymbol{X}(\varphi ,\theta )=\boldsymbol{X}_0+\int_{0}^{\varphi
}\left( r-a(s,\theta )\right) \left(
\begin{array}{c}
-\sin s \\
\cos s%
\end{array}%
\right) ds\cdot \left(
\begin{array}{c}
\boldsymbol{\Xi } \\
\boldsymbol{\Theta }(\theta)%
\end{array}%
\right) +h\left( \varphi ,\theta \right) \ \boldsymbol{\Psi }(\theta).
\label{simplef}
\end{equation}
We will also define%
\begin{equation}
\boldsymbol{X}_{oh}(\varphi ,\theta ):=\boldsymbol{X}_0+\int_{0}^{\varphi
}\left( r-a(s,\theta )\right) \left(
\begin{array}{c}
-\sin s \\
\cos s%
\end{array}%
\right) ds\cdot \left(
\begin{array}{c}
\boldsymbol{\Xi } \\
\boldsymbol{\Theta }(\theta )%
\end{array}%
\right) .  \label{xoh}
\end{equation}
Note that $\boldsymbol{U}(\varphi,\theta)=\cos\varphi~\boldsymbol{\Xi}+\sin\varphi~\boldsymbol{\Theta}(\theta)$
lies in the plane of $\boldsymbol{X}_{oh}(\mathsf{S})$ and is perpendicular to the cylinder
$\boldsymbol{X}_{oh}([0,2\pi],\theta)+\mathbb{R}~\boldsymbol{\Psi}(\theta)$. Since $G$ lies inside this cylinder
and $\boldsymbol{X}(\varphi,\theta)$ is a point of $\partial G$ on this cylinder, the vector
$\boldsymbol{U}(\varphi,\theta)$ is an outwards
normal to $\boldsymbol{X}(\mathsf{S})$ at $\boldsymbol{X}(\varphi,\theta)$.
In other words, writing $\boldsymbol{\tilde{X}}$ as in (\ref{XtildeX}) it follows that
\begin{equation*}
  \omega\mapsto\boldsymbol{\tilde{X}}(\omega):\mathbb{S}^2\to\partial G\subset\mathbb{R}^3
\end{equation*}
is the \lq inverse\rq\ of the Gauss map for $\partial G$ and hence
Lipschitz-continuous on $\mathbb{S}^2$ by Lemma \ref{techlip}.
Note that Lemma \ref{LLC} shows that Lipschitz-continuity of $\boldsymbol{\tilde{X}}$
on $\mathbb{S}^2$ implies Lipschitz-continuity of $\boldsymbol{X}$ on $\mathsf{S}$
(but not vice versa!). So we may state, allowing the notation $C^{0,1}_%
\mathrm{p}(\mathsf{S})$ as a restriction of $C^{0,1}_\mathrm{p}(\mathbb{R}^2)$, that
for each coordinate in (\ref{simplef}) the Lipschitz-continuity holds for:
\begin{eqnarray*}
\left( \varphi ,\theta \right) &\mapsto &\int_{0}^{\varphi }\left(
r-a(s,\theta )\right) \sin s\ ds\in C_{\mathrm{p}}^{0,1}\left( \mathsf{S}%
\right) , \\
\left( \varphi ,\theta \right) &\mapsto &\int_{0}^{\varphi }\left(
r-a(s,\theta )\right) \cos s~ds \cos \theta-h(\varphi,\theta)\sin\theta
\in C_{\mathrm{p}}^{0,1}\left( \mathsf{S}\right)  ,\\
\left( \varphi ,\theta \right) &\mapsto &\int_{0}^{\varphi }\left(
r-a(s,\theta )\right) \cos s~ds \sin \theta+h(\varphi,\theta)\cos\theta
\in C_{\mathrm{p}}^{0,1}\left( \mathsf{S}\right)  .
\end{eqnarray*}
Combining the second function above multiplied with $\sin\theta$ and the third multiplied with $\cos\theta$ we
find that $h$ is Lipschitz-continuous on $\mathsf{S}$, which is not sufficient for
Lipschitz-continuity on $\mathbb{S}^2$. The function
 $(\varphi,\theta)\mapsto(\sin\theta,\cos\theta)$ transferred to $\omega$ turns into a
 function that is not even continuous on $\mathbb{S}^2$.
 We need another proof that
\begin{equation}\label{hwig}
  \omega\mapsto\tilde{h}(\omega)\text{ for }\omega=\boldsymbol{U}(\varphi,\theta)
  \text{ and }\tilde{h}(\omega)=h(\varphi,\theta)
\end{equation}
is Lipschitz-continuous on $\mathbb{S}^2$ and for that we will use the next two lemmata. We continue this proof on page \pageref{cp}.
\end{proof}

Note that the sketch on the left of Fig.~\ref{3dbild} shows a domain with
boundary $\boldsymbol{X}_{oh}(\mathsf{S})$.\medskip

For the proof, that $\omega \mapsto \boldsymbol{\tilde{X}}_{oh}(\omega )$ from (\ref{xoh}) is
Lipschitz-continuous, we will use the following: if $\boldsymbol{\tilde X}(\mathbb{S}^2)$ describes the boundary of $G$, then
$\boldsymbol{\tilde X}_{oh}(\mathbb{S}^2)$ gives the boundary
of the 3d-shadow domain $Sh_{\boldsymbol{\Xi}}(G)$, when rotating around the central axis $\boldsymbol{\Xi}$.
The definition of 3d-shadow domain is found in Appendix \ref{Bp}.
The function $\boldsymbol{X}_{oh}$ does not only parametrize the
boundary of that 3d-shadow domain, but since
\begin{equation*}
  \boldsymbol{X}_{oh}(\varphi,\theta)=
  P_{\boldsymbol{\Psi}(\theta)}\left(\boldsymbol{X}_{oh}(\varphi,\theta)\right)=
  P_{\boldsymbol{\Psi}(\theta)}\left(\boldsymbol{X}(\varphi,\theta)\right)
\end{equation*}
we have that each shadow in the direction of $\boldsymbol{\Psi}(\theta)$ has the contour
parametrized by $\varphi\mapsto \boldsymbol{X}_{oh}(\varphi,\theta)$.

We start with an a-priori estimate for the position of $\partial G$ in relation with
the axes through the highest and lowest point of $G$.

\begin{lemma}
  Let $\boldsymbol{X}$ be as in (\ref{simplef}) with $\boldsymbol{X}_0$ as in (\ref{X0}).
  Then for each $(\varphi,\theta)\in\mathsf{S}$ one finds:
  \begin{equation}\label{rad}
\left|\boldsymbol{X}(\varphi,\theta)-\left\langle\boldsymbol{X}(\varphi,\theta),
\boldsymbol{\Xi}\right\rangle\boldsymbol{\Xi}\right|\le 2\boldsymbol{d}_G\,\left|\sin\varphi\right|.
  \end{equation}
\end{lemma}

\begin{proof}
In any horizontal direction the boundary $\partial G$ lies between the extreme cases of two-dimensional
curves of constant width. These extreme cases are the Reuleaux triangle
pointing left and the one pointing right. See Fig.~\ref{boundsfig}. Rotating the left image around the
vertical axis gives the area on the right, where $\partial G$ is located.

\begin{figure}[h]
\centering
\includegraphics[width=5cm]{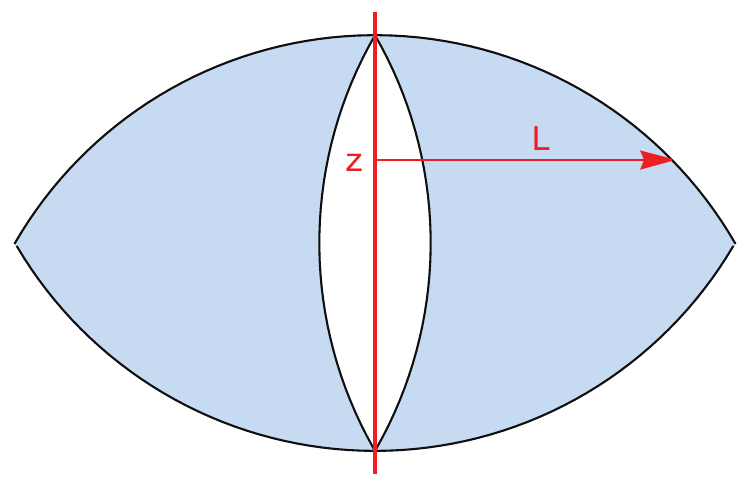}\hspace{1cm}\includegraphics[width=5.6cm]{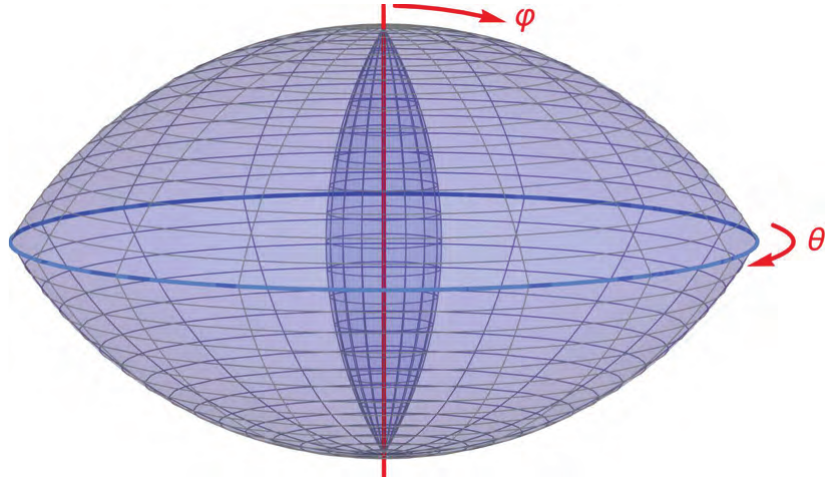}
\caption{The axis in red with the shaded parts showing the possible areas
for the projections in 2d. The shaded parts are the union of two Reuleaux
triangles minus their intersection. So with the assumption $r=1$ and $\left|z\right|\le r=1$
we find for the picture on the left  $L\in[-\sqrt{3},\sqrt{3}]$. If a body of constant width has $\left(0,0, \pm
1\right) $ on its boundaries, then the domain lies in
the rotated shaded part, which is sketched on the right.}
\label{boundsfig}
\end{figure}

So $\boldsymbol{X}(\varphi,\theta)$ is located
in the shaded area rotated around the $\boldsymbol{\Xi}$-axis. Setting
\begin{equation}\label{Z}
\boldsymbol{Z}(\varphi,\theta):=
\left\langle\boldsymbol{X}(\varphi,\theta),\boldsymbol{\Xi}\right\rangle\boldsymbol{\Xi}
\end{equation}
we find that for the top half, using $z\in[-1,1]$ for the relative height
$z:=\left(2/\boldsymbol{d}_G\right)\boldsymbol{Z}_3((\varphi,\theta))$ and the circle with center $(0,-1)^T$:
\begin{equation*}
\left|\boldsymbol{X}(\varphi,\theta)-\boldsymbol{Z}(\varphi,\theta)\right|\le
\boldsymbol{d}_G  \sqrt{4-(1+z)^2}
\end{equation*}
and for the bottom half, using the circle with center $(0,1)^T$:
\begin{equation*}
\left|\boldsymbol{X}(\varphi,\theta)-\boldsymbol{Z}(\varphi,\theta)\right|\le
\boldsymbol{d}_G  \sqrt{4-(1-z)^2}.
\end{equation*}
If $\left(\cos\hat{\varphi},\sin\hat{\varphi}\right)^T$ denotes the outward normal
direction for the curves
\begin{equation*}
  z\mapsto\left(z,\sqrt{4-(1+z)^2}\right) \text{ and }z\mapsto\left(z,\sqrt{4-(1-z)^2}\right),
\end{equation*}
then we find
\begin{equation}\label{XZ}
\left|\boldsymbol{X}(\varphi,\theta)-\boldsymbol{Z}(\varphi,\theta)\right|\le
2\boldsymbol{d}_G  \sin\hat{\varphi}.
\end{equation}
By the construction we find that $\boldsymbol{U}(\varphi,\theta)$ is an outward normal
in $\boldsymbol{X}(\varphi,\theta)$ for $\partial G$. Moreover, see (\ref{notp}), we have
$\boldsymbol{U}(\varphi,\theta)=\cos \varphi ~\boldsymbol{\Xi }+\sin \varphi ~\boldsymbol{\Theta }(\theta)$,
so $\boldsymbol{U}(\varphi,\theta)$ lies in the $\left\{\boldsymbol{\Theta}(\theta),\boldsymbol{\Xi}\right\}$-plane
that contains the curve $\varphi\mapsto\boldsymbol{X}_{oh}(\varphi,\theta)$ and is an outward normal to that curve.
A sketch of $\varphi\mapsto\boldsymbol{X}_{oh}(\varphi,\theta)$ in the
$\left\{\boldsymbol{\Theta}(\theta),\boldsymbol{\Xi}\right\}$-plane is found in Fig.~\ref{extreye}.

\begin{figure}[H]
  \centering
  \includegraphics[width=5.5cm]{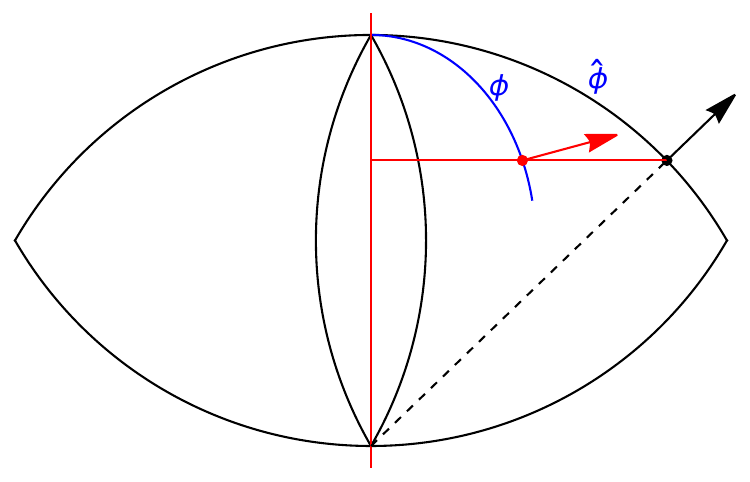}
  \caption{Comparing in the $\left\{\boldsymbol{\Theta}(\theta),\boldsymbol{\Xi}\right\}$-plane
   $\varphi\mapsto \boldsymbol{X}_{oh}(\varphi,\theta)$ for $\varphi\in[0,\frac12\pi]$
  and the circle parametrized by $\hat{\varphi}$.}\label{extreye}
\end{figure}

We have that if $\varphi\in(0,\tfrac12 \pi]$,  then $\hat{\varphi}\in (0,\varphi)$
and  if $ \varphi\in(\tfrac12 \pi,\pi)$,  then $\hat{\varphi}\in (\varphi,\pi)$.
Both cases imply that $\sin\hat{\varphi}\le\sin\varphi$. With similar estimates for $\varphi\in(\pi,2\pi)$
we may conclude from (\ref{XZ}) that
\begin{equation}\label{XZs}
\left|\boldsymbol{X}(\varphi,\theta)-\boldsymbol{Z}(\varphi,\theta)\right|\le
2\boldsymbol{d}_G  \left|\sin\varphi\right|,
\end{equation}
which is the claimed result.
\end{proof}

Note that $\left\langle\boldsymbol{X}(\varphi,\theta),\boldsymbol{\Xi}\right\rangle\boldsymbol{\Xi}
=\left\langle\boldsymbol{X}_{oh}(\varphi,\theta), \boldsymbol{\Xi}\right\rangle\boldsymbol{\Xi} $~
and  (\ref{rad}) implies that for all $(\varphi,\theta)\in\mathsf{S}$:
\begin{equation}\label{rado}
\left|\boldsymbol{X}_{oh}(\varphi,\theta)-\left\langle\boldsymbol{X}_{oh}(\varphi,\theta),
\boldsymbol{\Xi}\right\rangle\boldsymbol{\Xi}\right|\le \left|\boldsymbol{X}(\varphi,\theta)-\left\langle\boldsymbol{X}(\varphi,\theta),
\boldsymbol{\Xi}\right\rangle\boldsymbol{\Xi}\right|\le 2\boldsymbol{d}_G\,\left|\sin\varphi\right|.
  \end{equation}

\begin{lemma}
\label{shalem}Suppose that $\boldsymbol{X}$ from (\ref{simplef})
parametrizes the boundary of a body of constant width $G$ and
let $\boldsymbol{X}_{oh}$ be as in (%
\ref{xoh}). Then the function
\begin{equation}\label{xtiloh}
\omega \mapsto \boldsymbol{\tilde{X}}_{oh}(\omega ):\mathbb{S}%
^{2}\rightarrow \mathbb{R}^{3},
\end{equation}%
defined by $\boldsymbol{\tilde{X}}_{oh}(\boldsymbol{U}\left( \varphi ,\theta
\right) ):=\boldsymbol{X}_{oh}(\varphi ,\theta )$ for $(\varphi ,\theta )\in
\mathsf{S}$, is Lipschitz-continuous
and satisfies:

\begin{enumerate}
\item $P_{\boldsymbol{\Psi }(\theta )}\big(\boldsymbol{\tilde{X}}(\omega )%
\big)=\boldsymbol{\tilde{X}}_{oh}(\omega )$ for $\omega=\boldsymbol{U}( \varphi ,\theta)\in\mathbb{S}^2$,
and

\item $\partial \left(\mathit{Sh}_{\boldsymbol{\Xi }}(G )\right)=\boldsymbol{\tilde{X}}%
_{oh}(\mathbb{S}^{2})$.
\end{enumerate}
Here $P$ is as in (\ref{Pom}), $\boldsymbol{\Psi }(\theta),\boldsymbol{\Xi}$
as in (\ref{basis}),  $\boldsymbol{U}$ as in (\ref{defpar}) and
$\mathit{Sh}_{\boldsymbol{\Xi }}(\Omega )$ is the 3d-shadow
as in Definition \ref{shado3} and the rotation with respect to the axis \ $\boldsymbol{\Xi}$.
\end{lemma}

\begin{proof} We still assume (\ref{X0}).
From our construction one finds
that the function $\boldsymbol{\tilde{X}}_{oh}:\mathbb{S}^{2}\rightarrow
\mathbb{R}^{3}$ parametrizes the collection of boundaries of `2d-shadows' in
the directions $\boldsymbol{\Psi }(\theta )$ for
$\theta \in \left[ 0 ,\pi \right] $ and gives a bounded two-dimensional
manifold in $\mathbb{R}^{3}$. Each 2d-shadow
$P_{\boldsymbol{\Psi }(\theta )}(G)$ for $\theta \in \left[ 0,\pi \right] $
is a two-dimensional set of constant
width in the plane spanned by $\bf{\Xi}$ and $\bf{\Theta}(\theta)$. The 3d-domain $\Omega$
bounded by these curves, that is
\begin{equation*}
  \partial\Omega =\bigcup_{\theta \in \left[ 0 ,\pi \right] } P_{\boldsymbol{\Psi }%
(\theta )}\big(\boldsymbol{X}(\left[ 0,2\pi \right] ,\theta )\big)=\boldsymbol{\tilde{X}}_{oh}(\mathbb{S}^2),
\end{equation*}
is in general not a body of constant width and not even convex. But by rotating a body $G$ of
constant width around $\boldsymbol{\Xi}$,  we may use that each projection on
the $\left\{\boldsymbol{\Theta}(\theta),\boldsymbol{\Xi}\right\}$-plane
is a curve of constant width. For each fixed $\varphi$ the
function $\theta\mapsto \boldsymbol{X}_{oh}(\varphi,\theta)$
is Lipschitz-continuous according to Lemma \ref{drehlip} and with the estimate in
(\ref{rado}) we obtain:\begin{equation}\label{Lthet}
  \left|\boldsymbol{X}_{oh}(\varphi,\theta)-\boldsymbol{X}_{oh}(\varphi,\theta_0)\right|
  \le 2\boldsymbol{d}_G  \left|\theta-\theta_0\right|\left|\sin\varphi\right| \text{ for all }
  \theta,\theta_0 \text{ and }\varphi.
\end{equation}
Lipschitz-continuity of $\varphi\mapsto\boldsymbol{X}_{oh}(\varphi,\theta)$, with
constant $L=2\boldsymbol{d}_G$, follows from our 2d-construction:
\begin{equation}\label{Lvarp}
  \left|\boldsymbol{X}_{oh}(\varphi,\theta)-\boldsymbol{X}_{oh}(\varphi_0,\theta)\right|
  \le 2\boldsymbol{d}_G  \left|\varphi-\varphi_0\right| \text{ for all }
  \varphi,\varphi_0 \text{ and }\theta.
\end{equation}
For $\left|\theta-\theta_0\right|\le \frac12\pi$ we use a triangle inequality
with either $(\varphi,\theta_0)$ or $(\varphi_0,\theta)$ as an intermediate point and
both (\ref{Lthet}) and
(\ref{Lvarp}) to get  for all such $(\varphi,\theta),(\varphi_0,\theta_0)$:
\begin{equation}
  \left|\boldsymbol{X}_{oh}(\varphi,\theta)-\boldsymbol{X}_{oh}(\varphi_0,\theta_0)\right|
  \le 2\boldsymbol{d}_G \left( \left|\varphi-\varphi_0\right| +
  \left|\theta-\theta_0\right|\min\left(\left|\sin\varphi\right|,\left|\sin\varphi_0\right|\right)
  \right),\label{XohL}
\end{equation}
which fits with (\ref{Aeq1}) in Lemma \ref{LLC}.

For $\frac12\pi<\left|\theta-\theta_0\right|\le \pi$, and
$\varphi,\varphi_0$ both either near $0$ or $\pi$, one uses a triangle
inequality with $(0,\theta)$ or $(\pi,\theta)$, corresponding to the poles $\pm e_3$,
as an intermediate point and twice (\ref{Lvarp}).
Note that one may always choose $\theta,\theta_0$
such that one of these two cases holds, possibly by extending periodically as in
Definition \ref{funnydifdef} and one obtains estimates as in (\ref{Aeq2}).
So, with the equivalences in Lemma \ref{LLC}, the function  $\boldsymbol{\tilde{X}}_{oh}$
is Lipschitz-continuous on $\mathbb{S}^2$.
\end{proof}

\begin{proof}[Continued proof of Theorem \protect\ref{all}, a formula for $h$
and regularity]\label{cp}
Lemma \ref{techlip} states that $\omega\mapsto\boldsymbol{\tilde{X}}(\omega)$ is
Lipschitz-continuous  and Lemma  \ref{shalem} states Lipschitz-continuity for
$\omega\mapsto\boldsymbol{\tilde{X}}_{oh}(\omega)$.

To show the weighted Lipschitz continuity estimate as in (\ref{XohL}) for $h$, we use the expression
\begin{equation*}\label{hstart}
 h(\varphi,\theta)=\left( \boldsymbol{X}(\varphi,\theta)-
 \boldsymbol{X}_{oh}(\varphi,\theta)\right)\cdot\boldsymbol{\Psi}(\theta)
\end{equation*}
and the auxiliary term $\boldsymbol{Z}(\varphi,\theta)$ as in (\ref{Z}).
We find
\begin{gather*}
  \left|h(\varphi,\theta)-h(\varphi_0,\theta_0)\right|  \leq
  \left| \boldsymbol{X}(\varphi,\theta)-\boldsymbol{X}(\varphi_0,\theta_0)\right|
  \left|\boldsymbol{\Psi}(\theta)\right| +
  \left| \boldsymbol{X}_{oh}(\varphi,\theta)-\boldsymbol{X}_{oh}(\varphi_0,\theta_0)\right|
   \left|\boldsymbol{\Psi}(\theta)\right|~  \\
 \hspace{3cm}+~\left(\left| \boldsymbol{X}(\varphi_0,\theta_0)-\boldsymbol{Z}(\varphi_0,\theta_0)\right| +
   \left|\boldsymbol{X}_{oh}(\varphi_0,\theta_0)-\boldsymbol{Z}(\varphi_0,\theta_0)\right|\rule{0mm}{4mm}\right)
\left|\boldsymbol{\Psi}(\theta)-\boldsymbol{\Psi}(\theta_0)\right|\\
 \leq  \left| \boldsymbol{X}(\varphi,\theta)-\boldsymbol{X}(\varphi_0,\theta_0)\right|+
  \left| \boldsymbol{X}_{oh}(\varphi,\theta)-\boldsymbol{X}_{oh}(\varphi_0,\theta_0)\right| +
  4\, \boldsymbol{d}_G \left|\sin\varphi\right|\pi\left|\theta-\theta_0\right|.
\end{gather*}
In a similar way, we may show the estimate replacing $\sin\varphi$ by $\sin\varphi_0$
and hence with Lemma \ref{LLC} it follows that $\omega\mapsto\tilde{h}(\omega)$ as in (\ref{hwig})
is Lipschitz-continuous on $\mathbb{S}^2$.\medskip

Next we will derive the formula for $h$. Note that for $\boldsymbol{\Theta }$
and $\boldsymbol{\Psi }$ as functions of $\theta $:
\begin{equation*}
\boldsymbol{\Theta }(\theta +\varepsilon )=\cos \varepsilon \ \boldsymbol{%
\Theta }(\theta )+\sin \varepsilon \ \boldsymbol{\Psi }(\theta )\text{ \ and
\ }\boldsymbol{\Psi }(\theta +\varepsilon )=\cos \varepsilon \ \boldsymbol{%
\Psi }(\theta )-\sin \varepsilon \ \boldsymbol{\Theta }(\theta ).
\label{ThPsofeps}
\end{equation*}%
When there is no misunderstanding we skip the $\theta $-dependence of $%
\boldsymbol{\Theta }$\ and $\boldsymbol{\Psi }$ and use only $\boldsymbol{%
\Theta }=\boldsymbol{\Theta }(\theta )$ and $\boldsymbol{\Psi }=\boldsymbol{%
\Psi }(\theta )$ . Thus one computes
\begin{multline*}
\boldsymbol{X}(\varphi ,\theta +\varepsilon )-\boldsymbol{X}(\varphi ,\theta
)=\int_{0}^{\varphi }\left( r-a(s,\theta +\varepsilon )\right) \left(
\begin{array}{c}
-\sin s \\
\cos s%
\end{array}%
\right) ds\cdot \left(
\begin{array}{c}
\boldsymbol{\Xi } \\
\cos \varepsilon \ \boldsymbol{\Theta }+\sin \varepsilon \ \boldsymbol{\Psi }%
\end{array}%
\right) + \\
+h\left( \varphi ,\theta +\varepsilon \right) \left( \cos \varepsilon \
\boldsymbol{\Psi }-\sin \varepsilon \ \boldsymbol{\Theta }\right)
-\int_{0}^{\varphi }\left( r-a(s,\theta )\right) \left(
\begin{array}{c}
-\sin s \\
\cos s%
\end{array}%
\right) ds\cdot \left(
\begin{array}{c}
\boldsymbol{\Xi } \\
\boldsymbol{\Theta }%
\end{array}%
\right) -h\left( \varphi ,\theta \right) \ \boldsymbol{\Psi }
\\
=\int_{0}^{\varphi }\left( a(s,\theta )-a(s,\theta +\varepsilon )\right)
\left(
\begin{array}{c}
-\sin s \\
\cos s%
\end{array}%
\right) ds\cdot \left(
\begin{array}{c}
\boldsymbol{\Xi } \\
\boldsymbol{\Theta }%
\end{array}%
\right) +\left( h\left( \varphi ,\theta +\varepsilon \right) -h\left(
\varphi ,\theta \right) \right) \ \boldsymbol{\Psi \ }+ \\
-2\sin (\varepsilon /2)\left( h\left( \varphi ,\theta +\varepsilon \right)
\left(
\begin{array}{c}
\cos (\varepsilon /2) \\
\sin (\varepsilon /2)%
\end{array}%
\right) +\int_{0}^{\varphi }\left( r-a(s,\theta +\varepsilon )\right) \cos
sds\left(
\begin{array}{c}
\sin (\varepsilon /2) \\
-\cos (\varepsilon /2)%
\end{array}%
\right) \right) \cdot \left(
\begin{array}{c}
\boldsymbol{\Theta } \\
\boldsymbol{\Psi }%
\end{array}%
\right) .
\end{multline*}%
As $\boldsymbol{X}(\varphi ,\theta )$\ describes the surface of a body of
constant width $2r$ and
\begin{equation*}
\boldsymbol{X}(\varphi ,\theta )-\boldsymbol{X}(\varphi +\pi ,\theta )=2r%
\boldsymbol{U}(\varphi ,\theta )
\end{equation*}%
we find that for all $t\in\mathbb{R} $%
\begin{equation*}
\left\vert \boldsymbol{X}(\varphi ,\theta +t )-\boldsymbol{X}%
(\varphi +\pi ,\theta )\right\vert \leq 2r = \left\vert \boldsymbol{X}(\varphi
,\theta )-\boldsymbol{X}(\varphi +\pi ,\theta )\right\vert .  
\end{equation*}%
Note that%
\begin{equation*}
\left( \boldsymbol{X}(\varphi ,\theta +t )-\boldsymbol{X}(\varphi
+\pi ,\theta )\right) \cdot \boldsymbol{U}(\varphi ,\theta )=\left(
\boldsymbol{X}(\varphi ,\theta +t )-\boldsymbol{X}(\varphi ,\theta
)\right) \cdot \boldsymbol{U}(\varphi ,\theta )+2r
\end{equation*}%
and thus we necessarily have
\begin{equation}
\left( \boldsymbol{X}(\varphi ,\theta +t )-\boldsymbol{X}(\varphi
,\theta )\right) \cdot \boldsymbol{U}(\varphi ,\theta )\leq 0.  \label{cocol}
\end{equation}%
Since $\boldsymbol{U}(\varphi ,\theta )=\cos \varphi \ \boldsymbol{\Xi }%
+\sin \varphi \ \boldsymbol{\Theta }$ we find, using the
Lipschitz-continuity of $\theta \mapsto \boldsymbol{X}(\varphi ,\theta )$,
that%
\begin{gather*}
\left( \boldsymbol{X}(\varphi ,\theta +t )-\boldsymbol{X}(\varphi
,\theta )\right) \cdot \boldsymbol{U}(\varphi ,\theta )=\int_{0}^{\varphi
}\left( a(s,\theta )-a(s,\theta +t )\right) \sin \left( \varphi
-s\right) ds \\
-\sin \varphi \left( \sin t \ h\left( \varphi ,\theta +t
\right) +2\left( \sin (t /2)\right) ^{2}\int_{0}^{\varphi }\left(
r-a(s,\theta +t )\right) \cos s\ ds\right)
\end{gather*}%
\begin{equation*}
=t \left( \int_{0}^{\varphi }\frac{a(s,\theta )-a(s,\theta
+t )}{t }\sin \left( \varphi -s\right) ds-\sin \varphi \
h\left( \varphi ,\theta \right) \right) +\mathcal{O}\left( t
^{2}\right)\le 0 .
\end{equation*}%
For (\ref{cocol}) to hold it follows that for $\left|t\right| $ small:%
\begin{equation*}
\int_{0}^{\varphi }\frac{a(s,\theta )-a(s,\theta +t )}{t
}\sin \left( \varphi -s\right) ds-\sin \varphi \ h\left( \varphi ,\theta
\right) =\mathcal{O}\left( t \right) .
\end{equation*}%
And hence we find%
\begin{equation}
h\left( \varphi ,\theta \right) =\frac{\lim_{\varepsilon \rightarrow
0}\int_{0}^{\varphi }\frac{a(s,\theta )-a(s,\theta +\varepsilon )}{%
\varepsilon }\sin \left( \varphi -s\right) ds}{\sin \varphi }.\medskip
\label{coconode}
\end{equation}

It remains to show the regularity properties stated in the second item of
the theorem. These follow rather immediately. Whenever $a,\partial _{\theta
}a\in C_{\mathrm{p}}^{0}(\mathbb{R}^2)$ one finds from (\ref{coconode}) that%
\begin{equation*}
h\left( \varphi ,\theta \right) =-\frac{\int_{0}^{\varphi }a_{\theta }\left(
s,\theta \right) \sin \left( \varphi -s\right) ds}{\sin \varphi }
\end{equation*}%
as in (\ref{ha}). With $h$ satisfying (\ref{ha}) one finds for $a,\partial
_{\theta }a\in C_{\mathrm{p}}^{1}(\mathbb{R}^2)$ that also (\ref{hala}) is
satisfied.
\end{proof}

\appendix
\renewcommand{\thesection}{\Alph{section}} \setcounter{section}{0}


\section{On the distance in $\mathbb{S}^{2}$}

\label{Ap} Let $\boldsymbol{\tilde{f}}$ $:\mathbb{S}^{2}\rightarrow \mathbb{R%
}$ be some function. The standard definition for such a function $\boldsymbol{\tilde{f}}$
to be Lipschitz-continuous, is, that there exists $L>0$ such that
\begin{equation}
\left\vert \boldsymbol{\tilde{f}}(\omega )-\boldsymbol{\tilde{f}}(\omega
_{0})\right\vert \leq L\left\vert \omega -\omega _{0}\right\vert
\text{ for all }\omega ,\omega _{0}\in \mathbb{S}^{2}.  \label{LCZo}
\end{equation}
Since the functions we use are defined in terms of
$(\varphi,\theta)\in\mathsf{S}$ instead of $\omega\in\mathbb{S}^2$,
with $\mathsf{S}$ from (\ref{es}), we need to reformulate the
Lipschitz-condition in (\ref{LCZo}) to a condition for
\begin{equation*}
 \boldsymbol{f}:\mathsf{S}\rightarrow \mathbb{R} \text{ for }
 \omega =\boldsymbol{U}( \varphi ,\theta) \text{ and }
 \boldsymbol{f}(\varphi ,\theta )=\boldsymbol{\tilde{f}}(\omega )
\end{equation*}
with $\boldsymbol{U}$ from (\ref{defpar}). In other words, we have to replace
$\left|\omega-\omega_0\right|$ by an equivalent expression using
 $(\varphi,\theta)$ and $(\varphi_0,\theta_0)$. The corresponding estimates follow next.

\begin{lemma}
\label{LLC}Setting $\omega =\boldsymbol{U}\left( \varphi ,\theta \right) $
and $\omega _{0}=\boldsymbol{U}\left( \varphi _{0},\theta _{0}\right) $, one
finds that for all $(\varphi,\theta)$ and $(\varphi_0,\theta_0)$ in $\mathsf{%
S}$:

\begin{itemize}
\item if $\varphi ,\varphi _{0}\in \left[ 0,\pi \right] $ or $\varphi
,\varphi _{0}\in \left[ \pi ,2\pi \right] $:
\begin{equation}
\left\vert \omega -\omega _{0}\right\vert \leq \left\vert \varphi -\varphi
_{0}\right\vert +\left\vert \theta -\theta _{0}\right\vert \min \left(
\left\vert \sin \varphi \right\vert ,\left\vert \sin \varphi _{0}\right\vert
\right) \leq \pi \left\vert \omega -\omega _{0}\right\vert ;  \label{Aeq1}
\end{equation}

\item if $\varphi \in \left[ 0,\pi \right] $ and $\varphi _{0}\in \left[ \pi
,2\pi \right] $, or vice versa:%
\begin{equation}
\left\vert \omega -\omega _{0}\right\vert \leq \left\vert 2\pi -\varphi
-\varphi _{0}\right\vert +\left( \pi -\left\vert \theta -\theta
_{0}\right\vert \right) \min \left( \left\vert \sin \varphi \right\vert
,\left\vert \sin \varphi _{0}\right\vert \right) \leq \pi \left\vert \omega
-\omega _{0}\right\vert .  \label{Aeq2}
\end{equation}
\end{itemize}
\end{lemma}

\begin{figure}[H]
\centering
\includegraphics[width=.85\textwidth]{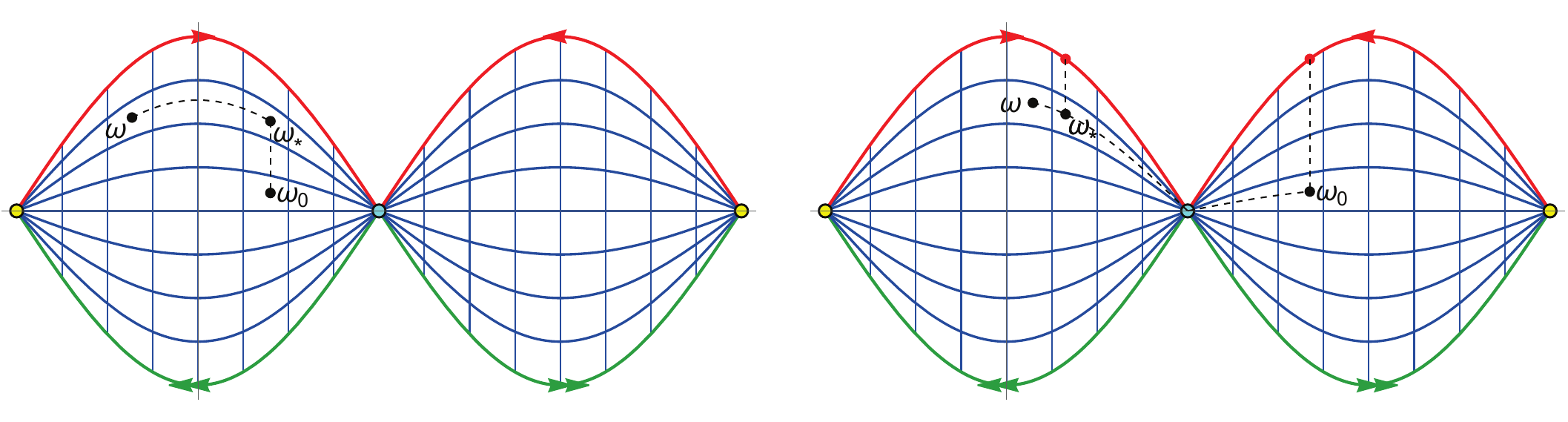}
\caption{A sketch of $\mathsf{S}$ deformed to show a distance equivalent to the one
in $\mathbb{S}^2$ identifying points at the upper boundary as well as points on the lower boundary.
For (\protect\ref{Aeq1}) see left and for (\protect\ref{Aeq2}) see
right. The variable $\protect\varphi\in [0,2\protect\pi ]$ moves from left
to right; $\protect\theta\left|\sin\protect\varphi\right|$ with $\protect%
\theta\in [-\frac12\protect\pi,\frac12 \protect\pi]$ increases in the
vertical direction.}
\label{connect}
\end{figure}

\begin{proof}
Assuming $\varphi ,\varphi _{0}\in \left[ 0,\pi \right] $ or $\varphi
,\varphi _{0}\in \left[ \pi ,2\pi \right] $ one considers as an intermediate
point $\omega _{\ast }=\boldsymbol{U}(\varphi _{0},\theta )$ and uses the
following estimates:

\begin{itemize}
\item The triangle inequality in $\mathbb{R}^3$: $\left\vert \omega -\omega
_{0}\right\vert \leq \left\vert \omega -\omega _{\ast }\right\vert
+\left\vert \omega _{\ast }-\omega _{0}\right\vert $.

\item Comparing the length via the circle with fixed $\varphi _{0}$ on the
sphere through the points $\omega $ and $\omega _{\ast }$ with the straight
line in $\mathbb{R}^{3}$ through those points gives:%
\begin{equation*}
\left\vert \omega -\omega _{\ast }\right\vert \leq \left\vert \varphi
-\varphi _{0}\right\vert \leq \tfrac{\pi }{2}\left\vert \omega -\omega
_{\ast }\right\vert .
\end{equation*}

\item A direct computation shows that%
\begin{equation*}
\left\vert \omega _{\ast }-\omega _{0}\right\vert =2\left\vert \sin \varphi
_{0}\right\vert \left\vert \sin \left( \tfrac{1}{2}\left( \theta -\theta
_{0}\right) \right) \right\vert
\end{equation*}%
and since $\theta -\theta _{0}\in \left[ -\pi ,\pi \right] $ one finds
\begin{equation*}
\tfrac{2}{\pi }\left\vert \theta -\theta _{0}\right\vert \leq 2\left\vert
\sin \left( \tfrac{1}{2}\left( \theta -\theta _{0}\right) \right)
\right\vert \leq \left\vert \theta -\theta _{0}\right\vert ,
\end{equation*}%
implying
\begin{equation}
\left\vert \omega _{\ast }-\omega _{0}\right\vert \leq \left\vert \theta
-\theta _{0}\right\vert \left\vert \sin \varphi _{0}\right\vert \leq \tfrac{%
\pi }{2}\left\vert \omega _{\ast }-\omega _{0}\right\vert .\label{refphi}
\end{equation}
By symmetry we may interchange $\omega$ and $\omega_0$ and hence
replace $\left\vert \sin \varphi _{0}\right\vert $ by $\left\vert \sin \varphi\right\vert $
and hence by  $\min\left( \left\vert \sin \varphi \right\vert ,\left\vert \sin
\varphi_{0}\right\vert \right) $ in (\ref{refphi}).

\item Both $\omega _{\ast }$ and $\omega _{0}$ lie on the circle on the unit
sphere with fixed $\varphi _{0}$. Since $\omega _{\ast }$ is the point on
that circle that is closest to $\omega $, one obtains%
\begin{equation*}
\left\vert \omega -\omega _{\ast }\right\vert \leq \left\vert \omega -\omega
_{0}\right\vert .
\end{equation*}%
A similar argument now for the circle on the unit sphere with fixed $\theta
_{0}$ shows%
\begin{equation}
\left\vert \omega _{\ast }-\omega _{0}\right\vert \leq \left\vert \omega
-\omega _{0}\right\vert .  \label{last}
\end{equation}
\end{itemize}

Combining these inequalities gives the estimates in (\ref{Aeq1}). \smallskip

For the second case we assume $\varphi \in \left[ 0,\pi \right] $, $\varphi
_{0}\in \left[ \pi ,2\pi \right] $ as in Fig.~\ref{connect} on the right. We
consider the shortest path from $\omega $ to $\omega _{0}$ through $\omega
_{\ast }=\boldsymbol{U}(2\pi -\varphi _{0},\theta )$ and the top or bottom
boundary of $\mathsf{S}$. Obviously $\left\vert \omega -\omega
_{0}\right\vert \leq \left\vert \omega -\omega _{\ast }\right\vert
+\left\vert \omega _{\ast }-\omega _{0}\right\vert $ still holds. As before
one finds%
\begin{equation*}
\left\vert \omega -\omega _{\ast }\right\vert \leq \left\vert 2\pi -\varphi
_{0}-\varphi \right\vert \leq \tfrac{\pi }{2}\left\vert \omega -\omega
_{\ast }\right\vert
\end{equation*}%
and since
\begin{equation*}
\left\vert \omega _{\ast }-\omega _{0}\right\vert =2\left\vert \sin \varphi
_{0}\right\vert \sin \left( \frac{\pi -\left\vert \theta -\theta
_{0}\right\vert }{2}\right)
\end{equation*}%
with $0\leq \pi -\left\vert \theta -\theta _{0}\right\vert \leq \pi $ one
obtains
\begin{equation*}
\left\vert \omega _{\ast }-\omega _{0}\right\vert \leq \left( \pi
-\left\vert \theta -\theta _{0}\right\vert \right) \left\vert \sin \varphi
_{0}\right\vert \leq \tfrac{\pi }{2}\left\vert \omega _{\ast }-\omega
_{0}\right\vert .
\end{equation*}%
Also as before we have%
\begin{equation*}
\left\vert \omega -\omega _{\ast }\right\vert \leq \left\vert \omega -\omega
_{0}\right\vert
\end{equation*}%
but the last inequality (\ref{last}) holds if $\left\vert \sin \varphi
_{0}\right\vert \leq \left\vert \sin \varphi \right\vert $. However, by taking
the minimum in (\ref{Aeq2}) the result holds true. The different cases are illustrated by
Fig.~\ref{connect}.
\end{proof}

\section{The inverse Gauss map for bodies of constant width}

As we mentioned in the introduction the Gauss map for bodies of constant width
$G$ is not necessarily uniquely defined on $\partial G$, but the inverse is.
This \lq inverse\rq\ is even Lipschitz and this can be found as a corollary in \cite{Ho}.
In the next Lemma we will give a short direct proof.

\begin{lemma}\label{techlip}
If $\boldsymbol{X}_1$ and $\boldsymbol{X}_2$ are two points on the
surface of a body $G\subset \mathbb{R}^{3}$
with constant width $\boldsymbol{d}_{G}$ and $\omega _{1}$, $\omega _{2}$
are outside normal directions at $\boldsymbol{X}_1$, $\boldsymbol{%
X}_2$, then%
\begin{equation}\label{b1}
\frac{\pi }{2}\left\vert \omega _{1}-\omega _{2}\right\vert \boldsymbol{d}%
_{G}\geq \left\vert \boldsymbol{X}_1-\boldsymbol{X}_2\right\vert .
\end{equation}
\end{lemma}

\begin{figure}[H]
  \centering
  \includegraphics[width=7cm]{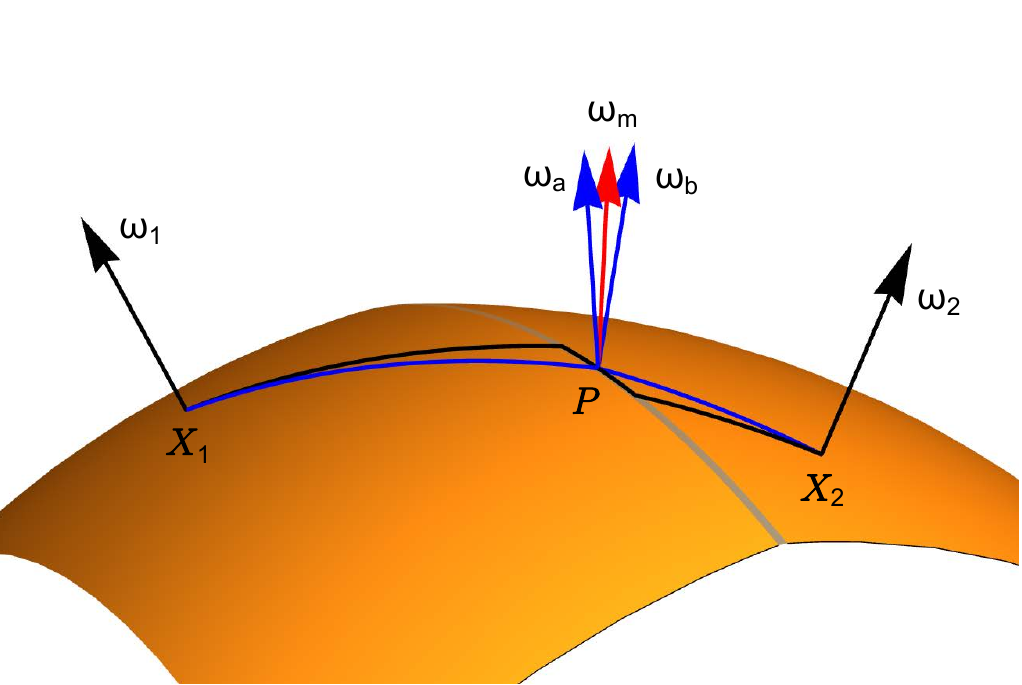}\hspace{1cm}\includegraphics[width=45mm]{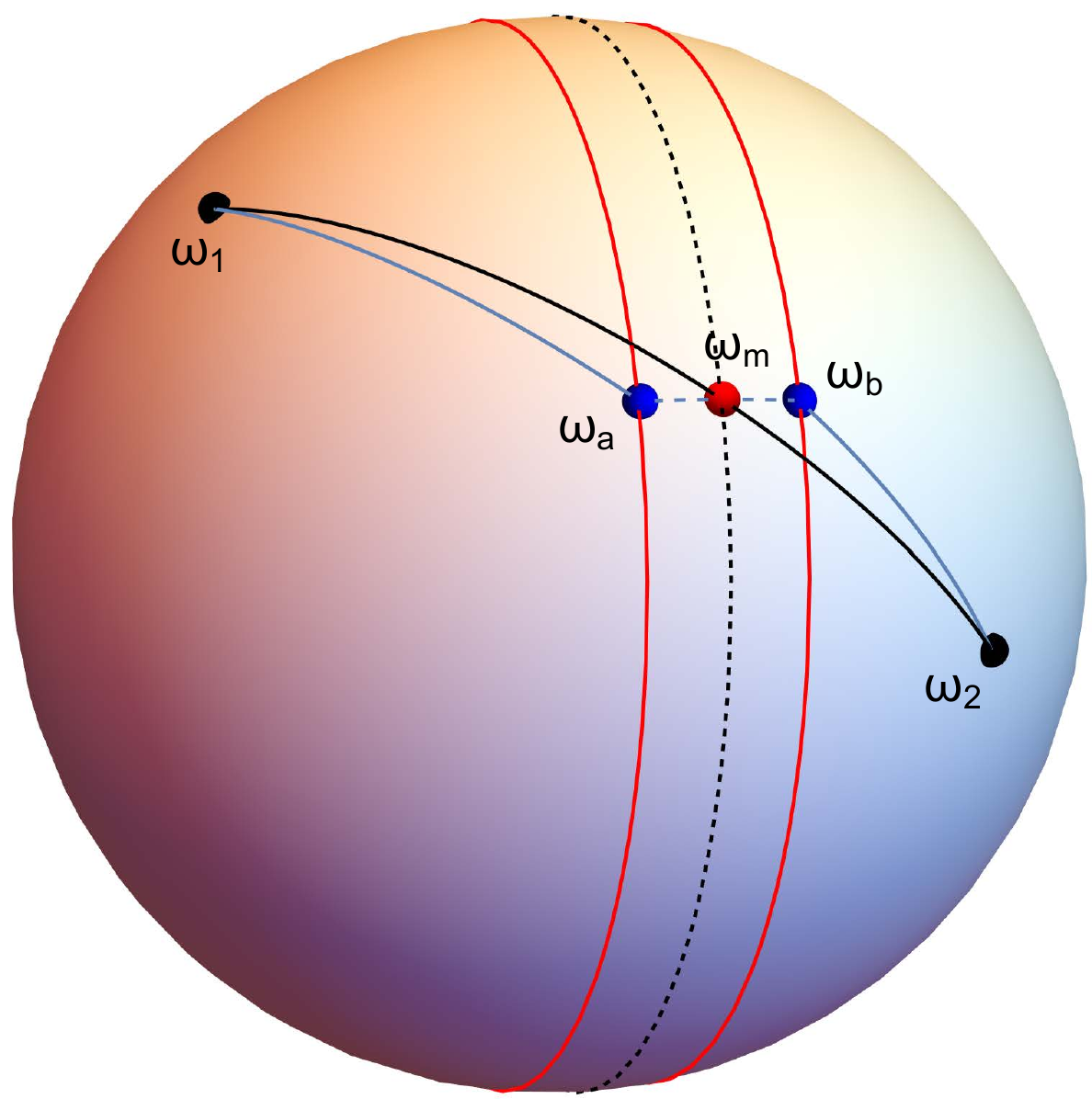}
  \caption{On the left parts of the two spheres from the proof of Lemma \ref{techlip}.
  $G$ lies below these two spheres and touches them in $\boldsymbol{X}_1$ and $\boldsymbol{X}_2$. On the right the
  various $\omega_{...}$ on $\mathbb{S}^2$ illustrating $\ell(\omega_1,\omega_2)=\ell(\omega_1,\omega_m)+\ell(\omega_m,\omega_2)\ge
  \ell(\omega_1,\omega_a)+\ell(\omega_b,\omega_2)$}\label{texlipfig}
\end{figure}

\begin{proof}
Let $\ell :\mathbb{S}^{2}\times \mathbb{S}^{2}\rightarrow \left[ 0,\pi %
\right] $ be the distance function on the sphere $\mathbb{S}^{2}$, that is%
\[
\ell \left( \omega _{\alpha},\omega _{\beta}\right) =\arccos \left( \omega _{\alpha}\cdot
\omega _{\beta}\right) \text{ for all }\omega _{\alpha},\omega _{\beta}\in \mathbb{S}^{2}.
\]%
Since $G$ is a body of constant width $\boldsymbol{d}_{G}$ it holds for any
$\boldsymbol{X}\in\partial G$ with outside normal $\omega$ that
\begin{equation*}
  x\in G \implies x\in \overline{B_{\boldsymbol{d}_G}}\left(\boldsymbol{X}-\boldsymbol{d}_G\omega\right).
\end{equation*}
Here $\overline{B_{r}}(\boldsymbol{M})$ is the closed ball of radius $r$ and
center $\boldsymbol{M}$. We call $\omega\in\mathbb{S}^2 $
an outside normal at $\boldsymbol{P}\in \partial G$ if $x\cdot \omega \leq
\boldsymbol{P}\cdot \omega $ for all $x\in G$.  Set
\begin{gather*}
  B_i:=B_{\boldsymbol{d}_{G}}\big(
  \boldsymbol{X}_i-\boldsymbol{d}_{G}\omega _{i}\big) \text{ for }i=1,2
,\\
L:=\overline{B_1}\cap\overline{B_2} \text{ and }C:=\partial B_1\cap\partial B_2.
\end{gather*}

Except for the circle $C$ where the two spheres intersect, there
is a unique outside normal direction on $\partial L$. There is a band on $\mathbb{S}%
^{2}$ that contains the outside normal directions connected to the circle $C$.
The circle in the middle of the band in $\mathbb{S}^2$ we call $S_{0}$. The
shortest path from $\omega _{1}$ to $\omega _{2}$ on $\mathbb{S}^{2}$
intersects $S_{0}$ at some $\omega _{m}$. There is a unique $\boldsymbol{P}\in C$
that connects with $\omega _{m}$. Set $\omega _{a}$ to be the outside normal
at $\boldsymbol{P}$ with respect to $\overline{B_1}$
and $\omega _{b}$ to be the outside normal at $\boldsymbol{P}$ with respect
to $\overline{B_2}$.

Since $\omega _{m}$ lies on the shortest path from $\omega _{1}$ to $\omega
_{2}$, we find that%
\[
\ell \left( \omega _{1},\omega _{2}\right) =\ell \left( \omega _{1},\omega
_{m}\right) +\ell \left( \omega _{m},\omega _{2}\right) .
\]%
Moreover $\ell \left( \omega _{1},\omega _{m}\right) \geq \ell \left( \omega
_{1},\omega _{a}\right) $ and $\ell \left( \omega _{m},\omega _{2}\right)
\geq \ell \left( \omega _{b},\omega _{2}\right) $. So one finds%
\begin{gather*}
\left\vert \boldsymbol{X}_1-\boldsymbol{X}_2\right\vert \leq
\left\vert \boldsymbol{X}_1-%
\boldsymbol{P}\right\vert +\left\vert \boldsymbol{P}-
\boldsymbol{X}_2\right\vert  \\
\leq \boldsymbol{d}_{G}\ \ell \left( \omega _{1},\omega _{a}\right) +%
\boldsymbol{d}_{G}\ \ell \left( \omega _{b},\omega _{2}\right) \leq
\boldsymbol{d}_{G}\ \ell \left( \omega _{1},\omega _{2}\right) \leq \frac{%
\pi }{2}\left\vert \omega _{1}-\omega _{2}\right\vert \boldsymbol{d}_{G},
\end{gather*}
which was the claim.
\end{proof}

\section{Shadow domains}

\label{Bp}

In order to show that a body of constant width has some minimal regularity
property, namely a kind of Lipschitz-continuity under rotation, we need a
geometrical argument. Such an argument follows from \lq observing the
shadows\rq\ during rotation. We did not find such a tool in the literature
and supply it here.

\begin{definition}
\label{shadowfun}Suppose that $\Omega \subset \mathbb{R}^{2}$ is a bounded,
simply connected domain with $0\in \Omega $. We define $R_{\Omega }:\mathbb{%
R\rightarrow R}^{+}$ by%
\begin{equation}
R_{\Omega }(\psi ):=\sup \left\{ x\cos \psi +y\sin \psi ;\ \binom{x}{y}\in
\Omega \right\}  \label{RO}
\end{equation}%
and the rotational shadow domain of $\Omega $ by%
\begin{equation*}
\mathit{Sh}(\Omega ):=\left\{ \binom{r\cos \psi }{r\sin \psi };0\leq
r<R_{\Omega }(\psi )\text{ and }\psi \in \left[ 0,2\pi \right] \right\} .
\end{equation*}
\end{definition}

\begin{figure}[h]
\centering\includegraphics[width=.9\textwidth]{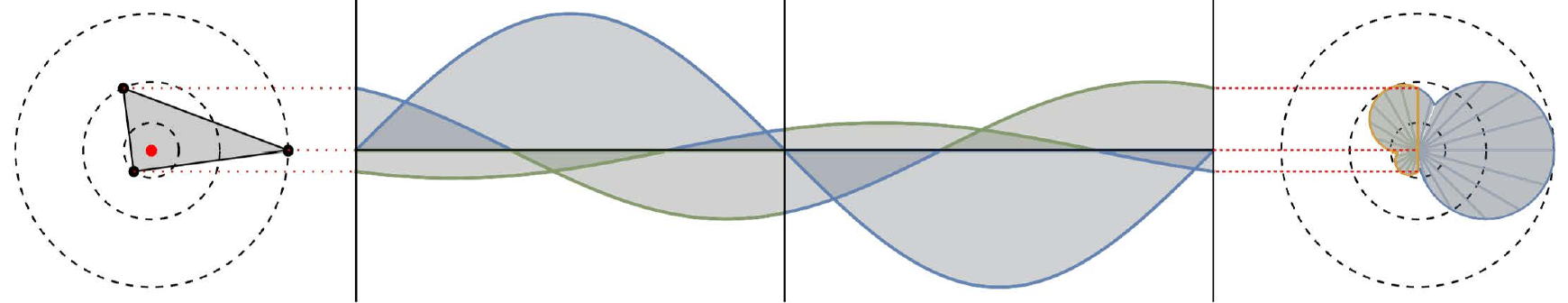}
\caption{On the left a triangle, in the middle $\protect\psi\mapsto R_{\Omega}(\protect\psi%
)$ as the maximum of the three functions related to the corners under counterclockwise rotation,
and on the right the rotational shadow domain of the triangle}
\label{shadowfig}
\end{figure}

The intersection of $\mathit{Sh}(\Omega )$ with the line
\begin{equation*}
\ell (\psi ):=\left\{ t\binom{\cos \psi }{\sin \psi };t\in \mathbb{R}\right\}
\end{equation*}%
gives precisely the shadow of $\Omega $ with the light at infinity in the
direction $\binom{-\sin \psi }{\cos \psi }$. See Fig.~\ref{shadowfig} in the
case of a triangle.

\begin{lemma}\label{drehlip}
Let $\Omega $ be as in Definition \ref{shadowfun}. The function $R_{\Omega }$
in (\ref{RO}) is Lipschitz-continuous with Lipschitz-constant at most%
\begin{equation}
L=\sup \left\{ \left\Vert x\right\Vert ;\ x\in \Omega \right\} .
\label{lippy}
\end{equation}
\end{lemma}

\begin{proof}
Let $co(\Omega )$ denote the convex hull of $\Omega $. It holds that $%
R_{co\left( \Omega \right) }(\psi )=R_{\Omega }(\psi )$. Note that taking
the convex hull also does not change $L$. Hence we may assume without loss
of generality that $\Omega $ is convex. The boundary of a bounded convex
domain in $\mathbb{R}^{2}$ with $0\in \Omega $ can be parametrized in polar
coordinates with $r(t)>0$ as follows:%
\begin{equation*}
\partial \Omega =\left\{ r(t)\binom{\cos t}{\sin t};t\in \left[ 0,2\pi %
\right] \right\} .
\end{equation*}%
For such a parametrization one finds%
\begin{align}
R_{\Omega }(\psi )& =\sup \left\{ r(t)\cos t\cos \psi +r(t)\sin t\sin \psi
;t \in \left[ 0,2\pi \right] \right\}  \notag \\
& =\sup \left\{ r(t)\cos \left( \psi -t\right) ;t \in \left[ 0,2\pi \right]
\right\} .  \label{ROOO}
\end{align}%
The function $\psi \mapsto r(t)\cos \left( \psi -t\right) $ is
Lipschitz-continuous with constant $\left\Vert r\right\Vert _{\infty }=L$ as
in (\ref{lippy}). A function defined as the supremum of Lipschitz-functions
with a uniform constant is Lipschitz-continuous with that same constant.
\end{proof}

Notice that (\ref{ROOO}) leads to
\begin{equation*}
R_{\Omega }(\psi )=\sup \left\{ r(\psi -s)\cos \left( s\right) ;\left\vert
s\right\vert <\tfrac{1}{2}\pi \right\} ,
\end{equation*}%
which again explains, why we call $\mathit{Sh}(\Omega )$ the rotational
shadow domain. Notice that, since $\cos s<0$ for $s\in [-\pi,-\frac{1}{2}\pi%
)\cup (\frac12\pi,\pi]$ and $0\in\Omega$, only the subinterval $%
(-\frac12\pi,\frac12\pi)$ contributes to this positive supremum.\medskip

One may extend this shadow in 2 dimensions to $3d$-shadows of a bounded
convex domain $\Omega \subset \mathbb{R}^{3}$. With the basis $\left\{
\boldsymbol{\Xi },\boldsymbol{\Theta }(\theta ),\boldsymbol{\Psi }(\theta
)\right\} $ as in (\ref{basis}) we define $P_{\boldsymbol{\Psi }(\theta )}:%
\mathbb{R}^{3}\rightarrow \mathbb{R}^{3}$, consistent with (\ref{Pom}), by%
\begin{equation*}
P_{\boldsymbol{\Psi }(\theta )}\left(
\begin{array}{c}
x_{1} \\
x_{2} \\
x_{3}%
\end{array}%
\right) :=\left\langle \boldsymbol{\Xi },x\rule{0mm}{3mm}\right\rangle
\boldsymbol{\Xi }+\left\langle \boldsymbol{\Theta }(\theta ),x\rule{0mm}{3mm}%
\right\rangle \boldsymbol{\Theta }(\theta )=\left(
\begin{array}{c}
\left( \cos \theta \ x_{1}+\sin \theta \ x_{2}\right) \cos \theta \\
\left( \cos \theta \ x_{1}+\sin \theta \ x_{2}\right) \sin \theta \\
x_{3}%
\end{array}%
\right) .
\end{equation*}
With the $\Xi$-axis being fixed the 3d-shadow is constructed as in the 2d-case
for each $\Xi$-coordinate being constant. One obtains a 3d-domain by joining the
2d-shadows from rotating around that axis. We use $\boldsymbol{\Xi}$ as before but the 3d-shadow
can be defined in any direction.

\begin{definition}
\label{shado3} Suppose that $\Omega \subset \mathbb{R}^{3}$ is strictly convex and bounded.
Then we define the $3d$-shadow domain in the directions perpendicular to the
$\boldsymbol{\Xi }$-axis by%
\begin{equation*}
\mathit{Sh}_{\boldsymbol{\Xi }}(\Omega ):=\bigcup \left\{ P_{\boldsymbol{%
\Psi }(\theta )}(\Omega );\left\vert \theta \right\vert \leq \tfrac{1}{2}\pi
\right\} .
\end{equation*}
\end{definition}

One may notice that this $3d$-shadow domain is related to the $2d$-shadows
for fixed $x_{3}$ through the formula
\begin{equation*}
\mathit{Sh}_{\boldsymbol{\Xi }}(\Omega )=\bigcup_{
x_{3}\in I}\left(
\begin{array}{c}
\mathit{Sh} \left(\left\{(x_1,x_2); (x_1,x_2,x_3)\in\Omega\right\}\rule{0mm}{3.5mm}\right)\\
x_3
\end{array}%
\right)   
\end{equation*}
for $I:=\left\{x_3; \exists (x_1,x_2,x_3)\in\Omega\right\}$.

\begin{figure}[p]
  \centering
  \includegraphics[width=\textwidth]{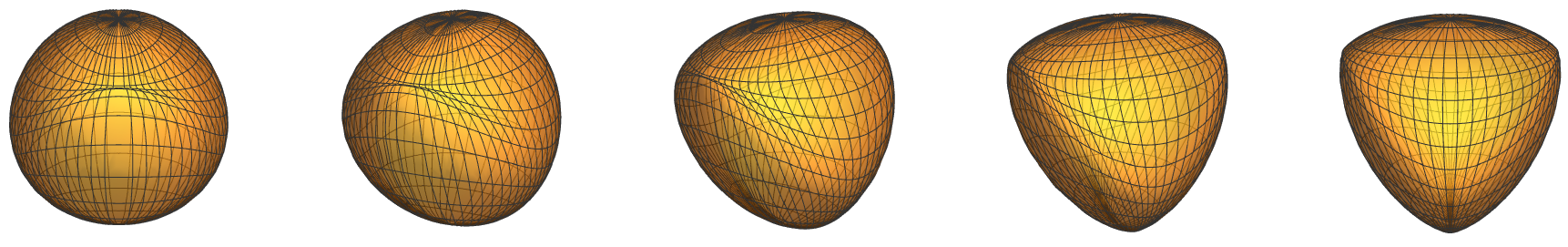}\\
  $a(\varphi,\theta)=-(\cos\theta)^2 \cos(3\varphi) $ with $r=r_0(a)=1$.\\[3mm]
  \includegraphics[width=\textwidth]{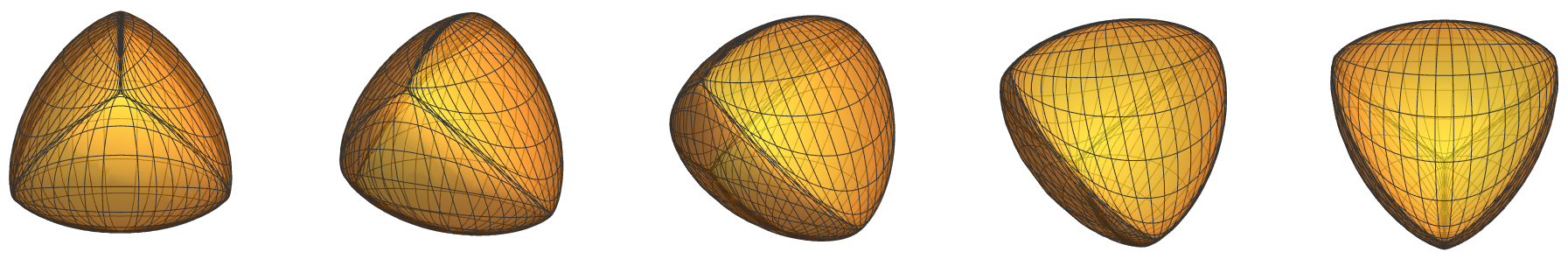}\\
  $a(\varphi,\theta)=-(\cos\theta)^2 \cos(3\varphi)+(\sin\theta)^2 \cos(3\varphi) $ with $r=r_0(a)=1.08867$.\\[3mm]
  \includegraphics[width=\textwidth]{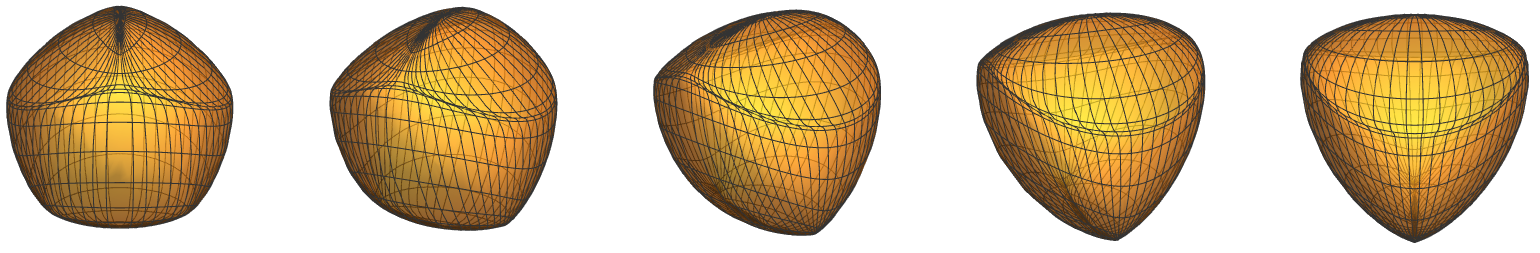}\\
  $a(\varphi,\theta)=-(\cos\theta)^2 \cos(3\varphi)+(\sin\theta)^2 \cos(5\varphi)  $ with $r=r_0(a)=1.11693$.\\[3mm]
  \includegraphics[width=\textwidth]{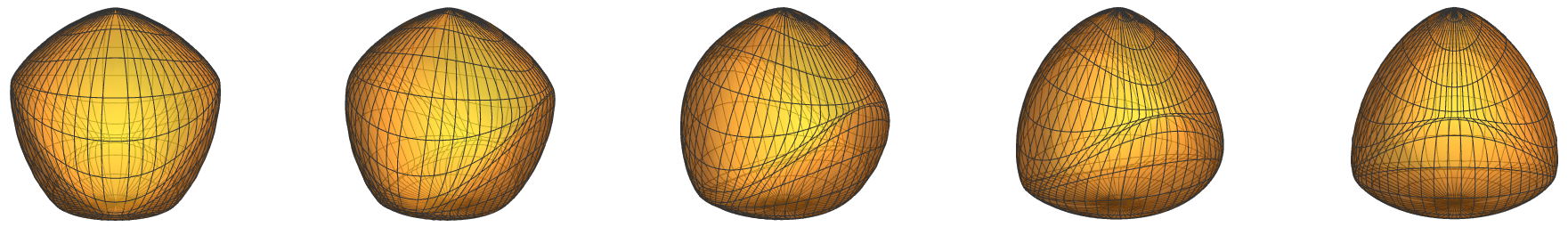}\\
  $a(\varphi,\theta)= (\cos\theta)^2 \cos(3\varphi)+(\sin\theta)^2 \cos(5\varphi)$ with $r=r_0(a)=1$.\\[3mm]
  \includegraphics[width=\textwidth]{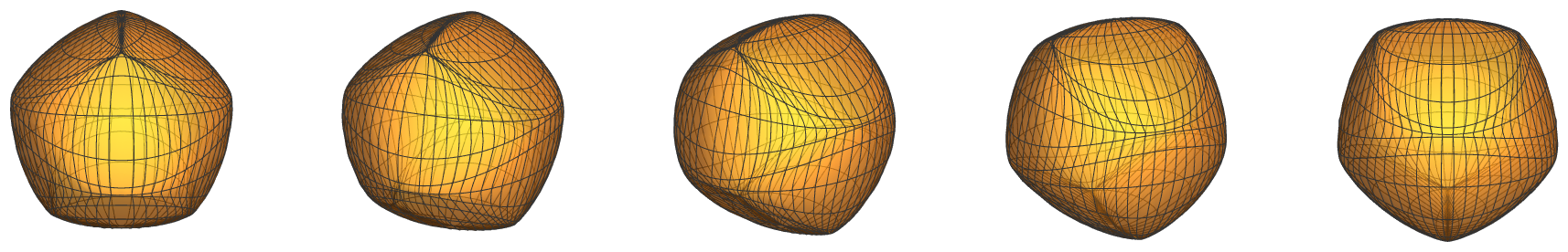}\\
  $a(\varphi,\theta)= -(\cos\theta)^2 \cos(5\varphi)+(\sin\theta)^2 \cos(5\varphi)$ with $r=r_0(a)= 1.01954$.\\[3mm]
  \includegraphics[width=\textwidth]{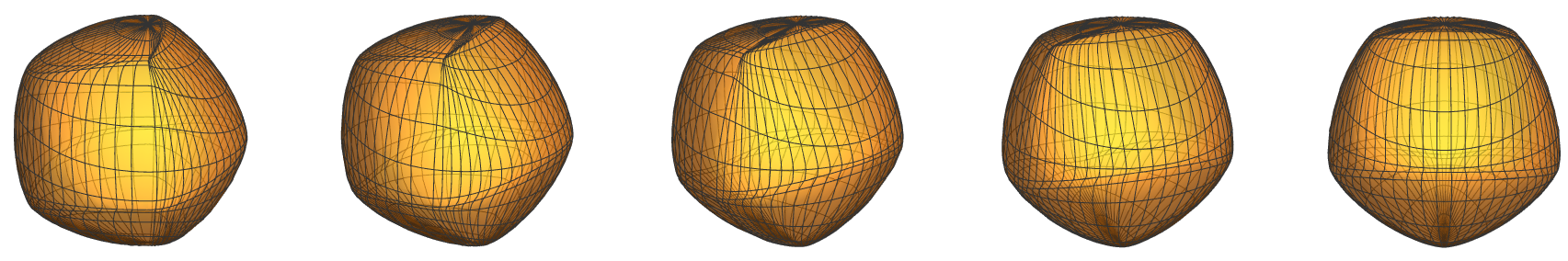}\\
  $a(\varphi,\theta)= -(\cos\theta)^2 \cos(5\varphi)+\left|\sin\theta\right|\sin\theta \sin(5\varphi) $ with $r=r_0(a)=1.1102$.
  \caption{Rotating views of 6 bodies of constant width; each connecting two simple curves of constant width}\label{exes}
\end{figure}

\end{document}